\documentclass[aap,MSNbibl,citesort,seceqn,secthm,dvips]{arximspdf}
\usepackage{graphics}
\doi{10.1214/09-AAP655}
\volume{20}
\issue{4}
\pubyear{2010}
\firstpage{1537}
\lastpage{1566}

\makeatletter

\newcommand{\eqref}[1]{(\ref{#1})}
\newtheorem {prop}{Proposition}[section]
\newtheorem {lemm}{Lemma}[section]
\newtheorem {cor}{Corollary}[section]
\newtheorem{assumption}{Assumption}

\newproclaim{rem}{Remark}[section]

\makeatother

\begin{document}
\begin{frontmatter}

\title{On the Wiener disorder problem}
\runtitle{On the Wiener disorder problem}

\begin{aug}
\author{\fnms{Semih Onur} \snm{Sezer}\ead[label=e1]{sezer@sabanciuniv.edu}\corref{}}
\runauthor{S. O. Sezer}
\affiliation{Sabanc{\i} University}
\address{Faculty of Engineering\\
\quad  and Natural Sciences \\
Sabanc{\i} University \\
Tuzla Istanbul 34956\\
Turkey \\
\printead{e1}} 
\end{aug}

\received{\smonth{4} \syear{2008}}
\revised{\smonth{8} \syear{2009}}

%
\begin{abstract}
In the Wiener disorder problem, the drift of a Wiener process
changes suddenly at some unknown and unobservable disorder time. The
objective is to detect this change as quickly as possible after it
happens. Earlier work on the Bayesian formulation of this problem
brings optimal (or asymptotically optimal) detection rules assuming
that the prior distribution of the change time is given at time
zero, and additional information is received by observing the Wiener
process only. Here, we consider a different information structure
where possible causes of this disorder are observed. More precisely,
we assume that we also observe an arrival/counting process
representing external shocks. The disorder happens because of these
shocks, and the change time coincides with one of the arrival times.
Such a formulation arises, for example, from detecting a change in
financial data caused by major financial events, or detecting
damages in structures caused by earthquakes. In this paper, we
formulate the problem in a Bayesian framework assuming that those
observable shocks form a Poisson process. We present an optimal
detection rule that minimizes a linear Bayes risk, which includes
the expected detection delay and the probability of early false
alarms. We also give the solution of the ``variational formulation''
where the objective is to minimize the detection delay over all
stopping rules for which the false alarm probability does not exceed
a given constant.
\end{abstract}

%
\begin{keyword}[class=AMS]
\kwd[Primary ]{62L10}
\kwd[; secondary ]{62L15}
\kwd{62C10}
\kwd{60G40}.
\end{keyword}

\begin{keyword}
\kwd{Sequential change detection}
\kwd{jump-diffusion processes}
\kwd{optimal stopping}.
\end{keyword}

\end{frontmatter}

\section{Introduction}
\label{secintro}
Suppose that at time $t = 0$ we start observing a Wiener process $X$ and
a simple Poisson process $N$ with arrival times
$(T_n )_{n \ge0}$.
The Poisson process
is assumed to apply
external shocks on $X$,
and these shocks will eventually cause a change in the drift of $X$.
The time $\Theta$, at which the drift changes is unknown and
unobservable. We only know that it
coincides with one of the arrival times according to the prior distribution
\begin{equation}
\label{distributionoftheta}
\hspace*{24pt} \mathbb{P}\lbrace\Theta= 0 \rbrace= \pi,\qquad
\mathbb{P}\lbrace\Theta= T_n \rbrace=(1- \pi)
(1-p)^{n-1}p \qquad
 \mbox{for all $n \ge1$}
\end{equation}
for some known $\pi\in[0,1)$ and $p \in(0,1]$. We also assume that
pre- and post-disorder drifts $\mu_0 $ and $\mu_1 $ are given, and the
arrival rate $\lambda$ of the Poisson process is known.

Our aim is to detect the time $\Theta$ as quickly as possible after it
happens, and by using our observations from the processes $X$ and $N$
only. 
More precisely, if we let $\mathbb{F}\equiv\{ \mathcal{F}_t \}_{t
\ge0}$ be the
observation filtration, our objective is to find
an $\mathbb{F}$-stopping time $\tau$ that minimizes the Bayes risk
\begin{equation}
\label{defR}
R(\tau, \pi) := \mathbb{P}\{ \tau< \Theta\} + c \cdot\mathbb
{E}(\tau- \Theta)^+
\end{equation}
for some delay cost $c>0$. If such a stopping time 
exists, then it resolves optimally the trade-off between early false
alarms and detection delay.

We also consider an alternative but related formulation, in which the
objective is to minimize the detection delay $ \mathbb{E}(\tau-
\Theta)^+$
over all $\mathbb{F}$-stopping times, for which the false alarm
frequency $\mathbb{P}
\{ \tau< \Theta\}$ is bounded above by 
a given constant $\alpha\in(0,1)$. Needless to say, this formulation
is more desirable if frequent false alarms cannot be tolerated.

Change detection problems have been studied in the literature with
numerous applications in different contexts. These applications include
quality control and fault detection in industrial processes, detection
of onset of an epidemic in biomedical signal processing, target
identification in national defense, intrusion detection in computer
networks and security systems,
threat detection in national security, pattern recognition in
seismology, detection of change in the riskiness of financial assets,
and many others. The reader may refer to 
\cite
{BassevilleNikiforovbook,Siegmundbook,Lai95,Lai01,Kent00,Shiryaevfinancialdata,WZS02,Baron02,TV04,TRBK06}, and the references therein for an extensive discussion
on these and other applications.\looseness=1

Earlier foundational studies on change detection problems include \cite{Lorden71} and \cite
{Page54} on non-Bayesian settings; and \cite
{GirshickRubin52} 
and \cite{Sh78}
on Bayesian formulations respectively. In particular, \cite{Sh78}
gives the
solution of the Bayesian 
formulation
of the Wiener disorder problem
for the Bayes risk in \eqref{defR} assuming that the change time has
an exponential prior distribution (see also \cite{Shiryaev63,Shiryaev65}).
Later, following \cite{Poor98}, this problem is reconsidered by \cite
{Beibel00}
for a different Bayes risk including an exponential penalty term (which
is more suitable for financial applications). Recently, \cite
{GapeevPeskir} obtained the solution of the finite-horizon version of
the original formulation of \cite{Sh78} (see also~\cite{PeskirShiryaevbook}, Chapter
6.22). The extension to the case where
observations consist of multiple Wiener processes is given by \cite
{dpsmultisource}. 

The common assumption in this line of work
is that
the change-time has (zero-modified) exponential distribution. Under
this assumption, the sufficient statistic (i.e., \textit{conditional
probability process}) is one dimensional, and it is possible to obtain
explicit results. In addition to this analytical advantage, the
exponential distribution can be regarded as a reasonable choice
for highly reliable systems considering the asymptotic approximation of
the exponential distribution with geometric distribution. That is, if
we perform independent experiments at times $\delta, 2 \delta, 3
\delta, \ldots,$  for $\delta> 0 $, where the failure (disorder)
probability is $\lambda\delta$, then as $\delta\to0^+$ we have
$\mathbb{P}(\mbox{time to first failure} > t) \to e^{- \lambda t}$.

In other settings where the prior distribution is not exponential, the
literature offers asymptotically optimal Bayes rules.
When the prior distribution is not exponential, sufficient statistics
are not one-dimensional anymore, and explicit results are difficult to
obtain, in which case asymptotically optimal rules prove useful for
online implementation.
The reader may refer to, for example,
\cite{Beibel94} and \cite{Beibel97}
for such asymptotical results including
explicit expansions of the optimal Bayes risks [which are modified
versions of \eqref{defR}];
see also \cite{PollakSiegmund85} for related results.
We refer the reader to the recent work \cite{BaronTartakovsky} for a
comprehensive asymptotical analysis of more general continuous-time
models (including the Wiener disorder problem).
The same work \cite{BaronTartakovsky} can also be consulted for a
brief survey and overview of the earlier work on
asymptotical detection theory.

In the aforementioned models,
the observed Wiener process is the only source of information for
detecting the change time. However, it is sometimes possible to observe
the external factors that are responsible for the disorder. This is
usually the case if
we would like
to detect, for example, a sudden change in financial data caused by
major financial events/news, or damages in structures caused by
earthquakes using continuously acquired vibration measurements (see
\cite{BassevilleNikiforovbook}, Chapters 1.2.5 and 11.1.4, for a
discussion on vibration monitoring in mechanical systems).
Here,
we consider such a setting where the underlying system is exposed to
observable shocks/impulses, and the disorder happens at one of these shocks.

\begin{figure}

\includegraphics{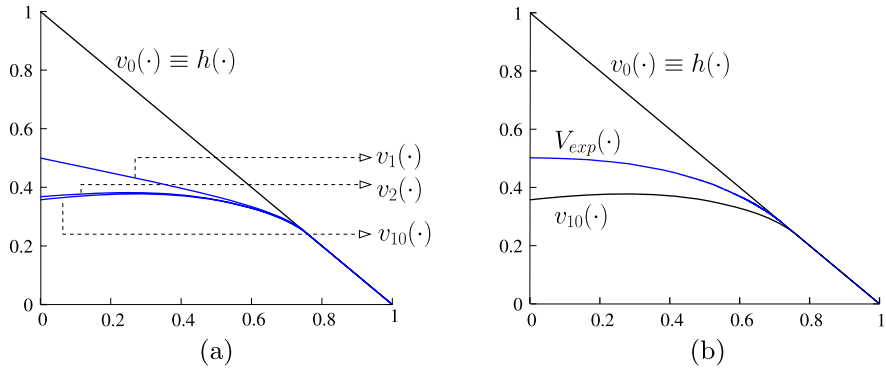}

\caption{We present a numerical example where $\mu=1$, $\lambda= 2$, $p =0.5$ and $c=0.5$.
Panel \textup{(a)} illustrates the sequential approximation of the optimal Bayes risk.
The functions $v_n (\cdot)$, for $n \le 10$, are computed by first finding the
threshold $\pi_n $ and then evaluating the exit time expectation $H_{\pi_n}[v_n] (\cdot)$ in
\protect\eqref{Hrexplicit} for $\pi \le \pi_n$.  The convergence is uniformly fast as given by
Lemma~\protect\ref{lemmuniformconvergence}. For $n=10$, we have $\| V - v_{10} \| \le 9.76 \cdot 10^{-4}$.
Panel \textup{(b)} compares two information levels on this detection problem. Recall that unconditional
distribution of the change time is exponential with parameter $\lambda p = 1$. If we only
observe the Wiener process (given this prior distribution) without observing the Poisson process $N$,
then we are in the framework considered by
  \protect\cite{PeskirShiryaevbook}, Section 6.22. The function $V_{\mathit{exp}}(\cdot)$ is the value function
  corresponding to this \textup{``}less information\textup{''} setting. It is computed by
  evaluating the expressions
  in \protect\cite{PeskirShiryaevbook}, pages 311--312, with the values of $\mu$ and $c$ given above. The figure
  illustrates that another observer who is also presented the process $N$ performs significantly
  better in detecting the change. }
    \label{figvs}
\end{figure}

Such a
formulation is considered for the first time by \cite{Moustakides08}
for a Brownian motion in a non-Bayesian framework, and under the
assumption that these shocks form a Poisson process. 
Sections 4 and 5 in \cite{Moustakides08} derive an optimal solution
for an (extended) Lorden criteria in terms of the (extended) CUSUM process.
However, to our knowledge, no Bayesian formulation of this problem has
been given yet. This formulation and its solution are the contributions
of the current paper.
It should be noted that under the distribution in \eqref
{distributionoftheta}, the unconditional distribution of $\Theta$ is
(zero-modified) exponential with parameter $\lambda p$. Hence, our
model can also be considered as a modification of the original
formulation in \cite{Sh78}. The major difference is that we not only
observe the underlying Wiener process but also the external causes
of the disorder. In this ``more informed'' setting, the detection
decision may improve greatly and this is indeed confirmed by our
numerical example in Figure~\ref{figvs}.

As an additional remark, we would like to note that although the change
can happen only at discrete points in time, a detection decision can be
made at any time. Hence, the problem is rather a continuous-time
problem as expected. It is essentially composed of a sequence of
hypothesis-testing problems: between two arrivals of the Poisson
process, the observer tests the hypotheses
\[
H_0\dvtx \mbox{drift } = \mu_0\quad    \mathrm{vs.}\quad   H_1\dvtx \mbox{drift } = \mu_1
\]
using the observations received from the Brownian motion. Indeed,
on every inter-arrival period $(T_n, T_{n+1})$, the conditional
probability process $\Pi_t := \mathbb{P}\{ \Theta\le t   | \mathcal
{F}_t \}$, for $t
\ge0$, follows the same dynamics as those of the sufficient statistic
$\widehat{\Pi}_t := \mathbb{P}\{ \mbox{$H_1$ is true}   | \mathcal
{F}_t \}$, for $t
\ge0$, of the sequential hypothesis-testing problem in \cite{Sh78}, Section
4.2; see Remark~\ref{remcomparingthedynamics}. If a decision
has not been made by the next arrival time $T_{n+1}$, then the
conditional probabilities are updated
and the hypothesis-testing problem restarts again with new (updated)
prior likelihoods.

In this paper, we show that the problem of minimizing the Bayes risk in
\eqref{defR} is equivalent to an optimal stopping problem in terms of
the conditional probability process 
$\Pi\equiv\{ \Pi_t \}_{t \ge0} $,
and it is optimal to stop the first time the process $\Pi$ exceeds a
threshold $\pi_{\infty}$.
The conditional probability process $\Pi$ is a jump-diffusion jointly
driven by the observed Wiener process and the Poisson process [see
\eqref{eqdynamicsofpi} for its dynamics]. To compute the optimal
threshold $\pi_{\infty}$ and the optimal Bayes risk, we transform the
corresponding optimal stopping problem into a sequence of stopping
problems for the diffusive part of the process $\Pi$. Each of these
sub-problems are solved by studying a free-boundary problem under a
smooth fit principle, and these solutions are then 
combined using a jump operator; see Sections \ref{secdpp} and \ref
{secsolution} below for details.
This approach is introduced for the first time by \cite
{dpsmultisource} in order to solve an optimal stopping problem
involving a discounted running cost only. 
In our setting, the problem includes a running cost and a terminal
cost, and involves no discounting. This requires nontrivial
modifications of their arguments as illustrated in Sections~\ref
{secdpp} and \ref{secsolution}.

In Section~\ref{secproblem} below, we formulate the problem
as an optimal stopping problem for the conditional probability process
$\Pi$, 
and we study the dynamics
of this process.
In Section~\ref{secdpp}, we introduce a jump operator whose role is to
incorporate the information generated by the Poisson process at every
arrival time. Using this operator, we construct the optimal Bayes risk
sequentially in Section~\ref{secsolution}, and we identify an optimal
Bayes rule. Finally, in Section~\ref{secvariational}, we solve the
variational formulation using the properties of the optimal solution
given in Section~\ref{secsolution}. Appendices at the end include some
of the lengthy derivations.

\section{Problem description}
\label{secproblem}
Let $(\Omega, \mathcal{H} , \mathbb{P})$ be a probability space hosting
a Wiener process $W$ and a simple Poisson process $N$ with arrival times
$(T_n )_{n \ge0}$
and rate $\lambda> 0$.
On this space, we have also an independent random variable $\zeta$
with the zero-modified geometric distribution
\begin{equation}
\label{distributionofzeta}
\quad \mathbb{P}\lbrace\zeta= 0 \rbrace= \pi,\qquad
\mathbb{P}\lbrace\zeta= n \rbrace=(1- \pi)
(1-p)^{n-1}p \qquad
 \mbox{for all $n \in\mathbb{N}$}
\end{equation}
for some $\pi\in[0,1)$ and $p \in(0,1]$.
In terms of these elements, we introduce a new $\mathbb{R}_+$-valued
variable 
\begin{equation}
\label{defTheta}
\Theta:= \sum_{i =0}^{\infty} T_i   1_{ \{ \zeta= i \} }
\end{equation}
representing the disorder time. Then, our observation process $X = \{
X_t \}_{t \ge0}$ can be defined as
\begin{equation}
\label{defX}
X_t := W_t + \mu\cdot(t - \Theta)^+\qquad \mbox{for all $t \ge0$.}
\end{equation}

In other words, as described in Section~\ref{secintro}, the process
$X$ is a Brownian motion gaining a drift $\mu$ at time $\Theta$, and
the change time $\Theta$ has zero-modified geometric distribution on
the Poissonian clock. With the notation in Section~\ref{secintro}, we
assume that $\mu_0 =0$ and $\mu_1 = \mu\ne0$ without loss of generality.

Let $\mathbb{F}\equiv \{ \mathcal{F}_t \}_{t \ge0}$ be the filtration of
the observed pair $(X,N)$;
that is, $\mathcal{F}_t := \sigma\{ X_s , N_s \dvtx s\le t \}$, for $t
\ge0$.
For an $\mathbb{F}$-stopping time $\tau$, let $R(\tau, \pi) $
denote the Bayes risk
\[
R(\tau, \pi) := \mathbb{P}^{\pi} \{ \tau< \Theta\} + c \cdot
\mathbb{E}^{\pi} (\tau
- \Theta)^+,
\]
in which
$\mathbb{P}^{\pi}$ is the probability measure $\mathbb{P}$ where
$\zeta$ has the distribution in \eqref{distributionofzeta}.
The Bayes risk above
includes the false alarm probability and the expected detection delay
cost for some $c > 0$. Our objective in this problem is to compute
\begin{equation}
\label{defU}
V( \pi) := \inf_{\tau\in\mathbb{F}} R(\tau, \pi),
\end{equation}
and if exists, find a stopping time attaining this infimum.

Using the standard arguments in \cite{Sh78}, Chapter 4, we can
transform the problem in \eqref{defU} into an optimal stopping problem
for the \textit{conditional probability process} defined as
\begin{equation}
\label{defPi}
\Pi_t := \mathbb{P}\{ \Theta\le t   | \mathcal{F}_t \}, \qquad   t
\ge0.
\end{equation}
More precisely, the minimal Bayes risk in \eqref{defU} is the value
function of the optimal stopping problem
\begin{equation}
\label{defV}
V( \pi) = \inf_{\tau\in\mathbb{F}}
\mathbb{E}^{\pi} \biggl[ \int_0^{\tau} g( \Pi_s ) \,  ds + h( \Pi
_{\tau} ) \biggr],
\end{equation}
where $g(\pi) := c \cdot\pi$ and $h(\pi) := 1-\pi$.

In Appendix~\hyperref[secappendixA]{A}, we show that the process $\Pi$ has the
characterization
\begin{eqnarray}
\label{eqPiandPhi}
\Pi_t &=& \frac{\Phi_t}{1 + \Phi_t },\qquad   \mbox{where}
\nonumber
\\[-8pt]
\\[-8pt]
\nonumber
\Phi_t :&=& \frac{L_t}{(1-p)^{N_t}} \Biggl( \frac{\pi}{1-\pi} +
\sum_{i=1}^{N_t} \frac{(1-p)^{i-1} p }{L_{T_i}} \Biggr)
\end{eqnarray}
in terms of 
\begin{equation}
\label{defL}
L_t := \exp\biggl\{ \mu X_t - \frac{\mu^2 }{2}t \biggr\}\qquad
\mbox{for $t \ge0$}.
\end{equation}
Using \eqref{eqPiandPhi} and \eqref{defL}, we obtain
\begin{equation}
\label{eqdynamicsofpi}
d \Pi_t = \mu\Pi_{t-} (1- \Pi_{t-}) \,d\widehat{W}_t + p (1- \Pi_{t-})\,dN_t,
\end{equation}
where
$
\widehat{W}_t := X_t - \mu\int_0^{t} \Pi_s \,ds,
$ 
for $t \ge0$, is a $(\mathbb{P}, \mathbb{F})$-Wiener process.
In Appendix~\hyperref[secappendixB]{B}, we also show that for $t \le s_1$ and $t \le s_2$ and
for $r, q \in\mathbb{R}$, we have
\begin{eqnarray}
\label{eqhatWandNindependent}
&&\mathbb{E}[ \exp\lbrace i   r ( \widehat{W}_{s_1}
- \widehat
{W}_{t} )
+ i   q ( N_{s_2} - N_{t} )
\rbrace| \mathcal{F}_t ]
\nonumber
\\[-8pt]
\\[-8pt]
\nonumber
&&\qquad= \exp\bigl\{
- \tfrac{1}{2} r^2 (s_1 - t) + \lambda(e^{iq} -1)(s_2 - t)\bigr\},
\end{eqnarray}
which implies that $\widehat{W}$ and $N$ are independent. This
further implies that the process~$\Pi$ whose dynamics are given in
\eqref{eqdynamicsofpi} is a strong Markov process.

\begin{rem} 
\label{remcomparingthedynamics}
Between two arrival times, the process $\Pi$ satisfies $d \Pi_t = \mu
\Pi_{t-} (1- \Pi_{t-})\, d\widehat{W}_t$, which coincides with the
dynamics of the conditional probability process in the sequential
hypothesis-testing problem considered in \cite{Sh78}, Section 4.2. In
that problem, an observer is given two hypotheses
\begin{equation}
\label{twohypotheses}
H_0 \dvtx \mbox{drift }= 0 \quad  \mbox{and} \quad   H_1 \dvtx \mbox{drift }=
\mu
\end{equation}
about an observed Wiener process. The hypotheses have
prior likelihoods $1-\pi$ and $\pi$ respectively, and the aim is to
identify the correct one as soon as possible.

In our problem, the change can happen only at one of the arrival times
of the Poisson process. Hence, between two arrival times $[T_n,
T_{n+1})$ the role of the process $\Pi$ is to indicate the posterior
likelihood of the hypothesis $H_1$, whose initial prior is $\Pi_{T_n}$
as of time $T_n$.
In this setting, if a decision is made by the next arrival time, the
hypothesis-testing problem terminates. Otherwise, it restarts with new
priors $1- \Pi_{T_{n+1}}$ and $\Pi_{T_{n+1}}$, respectively.
\end{rem}

\begin{rem}
\label{remlimitofpi}
Using its definition in \eqref{defPi}, it can easily be verified that
the process $\Pi$ is a bounded submartingale with a last element $\Pi
_{\infty} \le1$.
Moreover, thanks to bounded convergence theorem we have
\begin{eqnarray*}
1 &\ge&\mathbb{E}^{\pi} \Pi_{\infty} = \lim_{t \to\infty} \mathbb
{E}^{\pi} \Pi_t =
\lim_{t \to\infty} \mathbb{E}^{\pi} \bigl[ 1_{ \{ \Theta\le t \}
} \bigr] \\
&=&\lim_{t \to\infty} \mathbb{E}^{\pi} \bigl[ \mathbb{E}\bigl[ 1_{
\{ \Theta\le t \} }
| N_u : u\le t \bigr] \bigr]\\
& =&\lim_{t \to\infty} \mathbb{E}^{\pi} [ 1 - (1-\pi)
(1-p)^{N_t} ] =1 ,
\end{eqnarray*}
which implies that $\Pi_{\infty} =1 $, $\mathbb{P}^{\pi}$-a.s., for
all $\pi
\in[0,1]$.
\end{rem}

The limiting behavior of $\Pi$ implies that the exit time $\widetilde
{\tau}_r$ of $\Pi$
from an interval $[0,r)$, for $r \in[0,1)$, is finite $\mathbb
{P}^{\pi
}$-almost surely, for $\pi\in[0,1]$.
Indeed, the dynamics in \eqref{eqdynamicsofpi} give
\begin{eqnarray*}
&&\mathbb{E}^{\pi} \Pi_{t \wedge\widetilde{\tau}_r} = \pi+ \mathbb
{E}^{\pi} \biggl[
\int_0^{t \wedge\widetilde{\tau}_r } \mu\Pi_{u-} ( 1- \Pi_{u-} )\, d
\widehat{W}_u\\
&&\qquad\hspace*{68pt}{}+
\int_0^{t \wedge\widetilde{\tau}_r } p ( 1- \Pi_{u-} ) ( d N_u -
\lambda\, du )\\
 &&\qquad\hspace*{108pt}{}+ \int_0^{t \wedge\widetilde{\tau}_r } \lambda p ( 1-
\Pi
_{u-} )\, du
\biggr].
\end{eqnarray*}
Since $\Pi$ is bounded, the first two integral has zero expectations.
Moreover, for $u \le\widetilde{\tau}_r$, we have $1-\Pi_u \ge1- r $,
and this yields
\begin{equation}
\label{exittimeexpectationuniformlybounded}
1 \ge\mathbb{E}^{\pi} \Pi_{t \wedge\widetilde{\tau}_r} \ge
\lambda p ( 1- r ) \mathbb{E}^{\pi} t \wedge\widetilde{\tau}_r ,
\end{equation}
showing that $\mathbb{E}^{\pi} \widetilde{\tau}_r$ is uniformly
bounded, for
all $\pi\in[0,1]$, thanks to monotone convergence theorem.

\section{Dynamic programming operator}
\label{secdpp}
The first arrival time $T_1$ is a regeneration time of the conditional
probability process $\Pi$.
Therefore, if the process $\Pi$ has not been stopped yet, the minimal
Bayes risk that one can attain starting from $T_1$ should be computed
by evaluating the function $V(\cdot)$ at $\Pi_{T_1}$.
This tells us that the value function should satisfy the dynamic
programming equation
\begin{equation}
\label{DPPforV}
\hspace*{8pt}V(\pi) = \inf_{ \tau\in\mathbb{F}}
\mathbb{E}^{\pi} \biggl[ \int_0^{\tau\wedge T_1 } g( \Pi_t ) \,dt
+ 1_{ \{ \tau< T_1 \} } h( \Pi_{\tau} )
+ 1_{ \{ \tau\ge T_1 \} } V( \Pi_{T_1} ) \biggr].
\end{equation}

Until the first arrival time, $\Pi$ coincides with
a diffusion starting from $Y_0 = \pi$ and satisfying
\begin{equation}
\label{dynamicsofY}
d Y_t = \mu Y_{t} (1- Y_{t}) \,d\widehat{W}_t \qquad \mbox{for $t \ge0$}.
\end{equation}
Hence, a given stopping time $\tau\in\mathbb{F}$ should coincide on
the event
$\{ \tau< T_1\}$ with another stopping time of the process
$Y$.
This observation 
suggests that the function~$V(\cdot)$ should be a fixed point of the operator
%
\begin{eqnarray}
\label{defJ}
&&J[w](\pi) := \inf_{ \tau\in\mathbb{F}^{Y} } \mathbb{E}^{\pi}
\biggl[ \int_0^{\tau
\wedge T_1 } g( Y_t )\,   dt
+ 1_{ \{ \tau< T_1 \} } h( Y_{\tau} )
\nonumber
\\[-8pt]
\\[-8pt]
\nonumber
&&\qquad\hspace*{72pt}{}+ 1_{ \{ \tau\ge T_1 \} } w\bigl( Y_{T_1} + p (1-Y_{T_1}) \bigr)
\biggr],
\end{eqnarray}
which is obtained by replacing $\mathbb{F}$ in \eqref{DPPforV} with the
filtration $\mathbb{F}^Y$ of the process $Y$, and
$V(\cdot)$ with a bounded function $w(\cdot)$ on $[0,1]$.
Using the independence of $\widehat{W}$ and $N$, and the distribution
of $T_1$ we can rewrite this operator
\begin{equation}
\label{Jreloaded}
J[w](\pi) = \inf_{ \tau\in\mathbb{F}^Y}
\mathbb{E}^{\pi} \biggl[ \int_0^{\tau} e^{-\lambda t } [ g(
Y_t ) + \lambda
w( \mathbb{S}(Y_t) ) ]\, dt + e^{ - \lambda\tau} h( Y_{\tau} )
\biggr],
\end{equation}
where $\mathbb{S}(\pi) := \pi+ p (1-\pi)$.

In this section, we study the properties of the operator $J$ for a
suitable class of function $w(\cdot)$'s. Under certain assumptions on
$w(\cdot)$, we show that the infimum in~\eqref{Jreloaded} is attained
by the exit time of the process $Y$ from an interval of the form
$[0,r)$, and that the function $J[w](\cdot)$ solves the variational
inequalities of the optimal stopping problem in \eqref{Jreloaded}.
Using the results of this section, we show in Section~\ref
{secsolution} that the function $V(\cdot)$ indeed satisfies $V(\cdot
)=J[V](\cdot)$ as expected.\looseness=1

\begin{rem}
\label{rempropertiesofJ}
For a bounded function $w\dvtx [0,1] \mapsto\mathbb{R}_+$, we have $0 \le
J[w](\cdot) \le h(\cdot)$. Moreover for two bounded functions
$w_1(\cdot
) \le w_2(\cdot)$, we have $J[w_1](\cdot) \le J[w_2](\cdot)$.
\end{rem}
\begin{pf}
The upper bound $J[w](\cdot) \le h(\cdot)$ follows by taking
$\tau= 0$ in \eqref{Jreloaded}. Nonnegativity of $J[w](\cdot)$ and
the monotonicity of $w \mapsto J[w]$ are obvious.
\end{pf}

\subsection[Solving the optimal stopping problem in (3.4)]{Solving the optimal stopping
problem in \protect\eqref{Jreloaded}}
Below we solve the minimization problem in \eqref{Jreloaded} under the
following assumption.

\renewcommand{\theassumption}{A\arabic{assumption}}
\setcounter{assumption}{0}
\begin{assumption}\label{assa1} The function $w(\cdot)$ is an
arbitrary (but fixed) nonnegative and continuous function on $[0,1]$
bounded above by 
$h(\cdot)$.
\end{assumption}

Let us define the functions
\begin{eqnarray}
\label{defpsiandeta}
\psi(\pi) &:=& \pi^{m_1} (1-\pi)^{1-m_1}\quad   \mbox{and}
\nonumber
\\[-8pt]
\\[-8pt]
\nonumber
\eta(\pi) &:=& \pi^{m_2} (1-\pi)^{1-m_2}\qquad  \mbox{for $\pi\in
[0,1]$},
\end{eqnarray}
where $m_1 > 1$ and $m_2 < 0$ are the roots of the quadratic equation
\begin{equation}
\label{quadraticequationform}
m (m-1) = \frac{2 \lambda}{\mu^2}.
\end{equation}
The functions $\psi(\cdot)$ and $\eta(\cdot)$ are respectively, the
increasing and decreasing solutions (up to multiplication by a
constant) of the equation
$\mathcal{A}_0 f(\pi) = \lambda f(\pi)$,
where $\mathcal{A}_0 $ is the infinitesimal generator of the diffusion process
$Y$ in \eqref{dynamicsofY}; that is, $\mathcal{A}_0 f(\pi) :=
\frac{1}{2} \mu
^2 \pi^2 (1-\pi)^2 f''(\pi)$. It is easy to verify that the functions
$\psi(\cdot)$ and $\eta(\cdot)$ satisfy the boundary conditions
\begin{eqnarray}
\label{boundarybehaviorofpsiandeta}
\psi(0+) &=& \psi'(0+) = 0 = \eta'(1-)= \eta(1-),
\nonumber
\\[-8pt]
\\[-8pt]
\nonumber
\psi(1-) &=& \psi'(1-) = \infty= \eta'(0+)= \eta(0+) ,
\end{eqnarray}
and that their Wronskian is $m_1 - m_2$.

In terms of the drift $r(\cdot) \equiv0$ and the volatility $\sigma
(\pi) = \mu\pi(1-\pi)$ of
the process~$Y$ [see the dynamics in \eqref{dynamicsofY}], let
$S(\cdot)$ denote the \textit{scale function}
\[
S(\pi) := \int _{d}^{\pi} S(dy),\qquad  \mbox{where }
S(dy):=\exp \biggl\{ \int _c^{y} \frac{r(z) }{ \sigma^2(z)}\,dz
\biggr\} \,dy = dy
\]
for arbitrary $c,d \in(0,1)$, and let $M(\cdot)$ be the \textit{speed measure}
\[
M(dy) := \frac{dy}{\sigma^2(y) S'(y)} = \frac{dy}{\mu^2 \pi^2
(1-\pi)^2}
\]
for $y \in(0,1)$.
Feller boundary test at the right boundary $\{1\}$ gives
\begin{equation}
\label{Fellertestat1} \int _c^1 \int _y^1 S(dz) M(dy)
= \infty, \qquad \int _c^1 \int _y^1 M(dz) S(dy) =
\infty,
\end{equation}
and according to \cite{KarlinTaylorbook}, Table 6.2, we conclude that
the right boundary is \textit{natural}.
This implies that the process $Y$ cannot reach the right boundary in
finite time.
On the other hand,
the process $\{ 1-Y_t \}_{t \ge0}$ has the same dynamics in \eqref
{dynamicsofY} and by symmetry 
the left boundary $\{0\}$ is also natural for $Y$. Indeed, by a change
of variable in \eqref{Fellertestat1} as $u=1-y$ and $q = 1-z $ we
get the Feller boundary test at~$\{0\}$\looseness=1
\begin{eqnarray*}
\infty&=& \int _c^1 \int _y^1 S(dz) M(dy) =
\int _{0}^{1-c} \int _0^u dq \frac{du}{\mu^2 u^2
(1-u)^2}\\
& =& \int _{0}^{1-c} \int _0^u S(dq) M(du),
\\
\infty&=& \int _c^1 \int _y^1 M(dz) S(dy) =
\int _0^{1-c} \int _0^u \frac{dq}{\mu^2 q^2 (1-q)^2} \,du\\
&=& \int _0^{1-c} \int _0^u M(dq) S(du) ,
\end{eqnarray*}
which gives the same conclusion for the left boundary.

\begin{rem}
\label{remintegrabilitycond}
The process $Y$ is a bounded martingale
[see \eqref{dynamicsofY}],
and we have
\begin{eqnarray*}
\mathbb{E}^{\pi} \biggl[ \int_0^{\infty} e^{-\lambda t } [ g(
Y_t ) +
\lambda w( \mathbb{S}(Y_t) ) ]\, dt \biggr] &\le&\mathbb{E}^{\pi
} \biggl[ \int
_0^{\infty} e^{-\lambda t } [ g( Y_t ) + \lambda\| w \| ]
\,dt \biggr] \\
&=& \int_0^{\infty} e^{-\lambda t } [ \mathbb{E}^{\pi} g( Y_t )
+ \lambda\|
w \| ] \,dt\\
&=& \int_0^{\infty} e^{-\lambda t } [ g( \pi) + \lambda\| w \|
]\, dt < \infty.
\end{eqnarray*}
\end{rem}

\begin{lemm}
\label{lemmexittimeexpectation1}
For $0< l \le r < 1$, and let $\tau_{l,r}$ be the exit time of the
process $Y$ from the interval $(l,r)$. The expectation
\begin{equation}
\label{defHlr}
\qquad H_{l,r}[w](\pi) := \mathbb{E}^{\pi} \biggl[ \int_0^{\tau_{l,r}}
e^{-\lambda t }
[ g( Y_t ) + \lambda w( \mathbb{S}(Y_t) ) ] \,dt + e^{ -
\lambda\tau
_{l,r} } h( Y_{\tau_{l,r}} ) \biggr],
\end{equation}
has the explicit form
\begin{eqnarray}
\label{Hlrexplicit}
 H_{l,r}[w](\pi) &=& \psi(\pi)\biggl[ C_1 + \int_{\pi}^r u_1[w](y)\,dy
\biggr]\nonumber\\
&&{}+ \eta(\pi)\biggl[ C_2 - \int_{\pi}^r u_2[w](y)\,dy \biggr]\\
&&{}+ h(l) \frac{\psi(\pi) \eta(r) - \psi(r) \eta(\pi) }{\psi
(l)\eta(r)-
\psi(r)\eta(l)}
+ h(r) \frac{\psi(l) \eta(\pi) - \psi(\pi) \eta(l) }{\psi
(l)\eta(r)-
\psi(r)\eta(l)}\nonumber
\end{eqnarray}
for $\pi\in(l,r)$, in terms of
\begin{eqnarray}
\label{defuandC}
u_1[w] (y) &:=& 2 \frac{g(y) + \lambda w(\mathbb
{S}(y))}{(m_1-m_2)\sigma^2(y)}
\eta(y),
\nonumber\\
u_2[w] (y) &:=& 2 \frac{g(y) + \lambda w(\mathbb{S}(y))}{(m_1-m_2)
\sigma^2(y)} \psi(y),
\nonumber
\\[-8pt]
\\[-8pt]
\nonumber
C_2 &:=& \frac{\eta(l) \int_l^r u_2[w](y)\, dy   - \int_l^r u_2[w](y)
\,dy\, \psi(l) }{\psi(r)\eta(l)- \psi(l)\eta(r)}   \psi(r),
\\ C_1 &:=& - \frac{\eta(r)}{\psi(r)} C_2.\nonumber
\end{eqnarray}
%
%
%
Clearly, $H_{l,r}[w](\cdot)$ is nonnegative, and we have
$H_{l,r}[w](\pi
)=h(\pi)$, for $\pi\notin(l,r)$.
\end{lemm}
%

%
\begin{pf}
Nonnegativity of $H_{l,r}[w](\cdot)$ and the identity
$H_{l,r}[w](\cdot
)=h(\cdot)$, on $[0,1] \setminus(l,r)$, are obvious. For $\pi\in
(l,r)$, let $f(\cdot)$ denote the function on the right-hand side in
\eqref{Hlrexplicit}. Direct computation shows that $f(\cdot)$ satisfies
\[
( - \lambda+ \mathcal{A}_0 )f(\pi) + g(\pi) + \lambda w(\mathbb
{S}(y)) =0\qquad
\mbox{on $\pi\in(l,r)$},
\]
with boundary conditions $f(l+)= h(l) $ and $f(r-)= h(r) $. Moreover,
its derivative (with respect to $\pi$) is
\begin{eqnarray*}
&&\psi'(\pi)\biggl[ C_1 + \int_{\pi}^r u_1[w](y)\,dy \biggr]
+ \eta'(\pi)\biggl[ C_2 - \int_{\pi}^r u_2[w](y)\,dy \biggr]\\
&&\qquad{}+ h(l) \frac{\psi'(\pi) \eta(r) - \psi(r) \eta'(\pi) }{\psi
(l)\eta(r)-
\psi(r)\eta(l)}
+ h(r) \frac{\psi(l) \eta'(\pi) - \psi'(\pi) \eta(l) }{\psi
(l)\eta(r)-
\psi(r)\eta(l)},
\end{eqnarray*}
which is bounded on $[l,r]$. Also, observe that the exit time $\tau
_{l,r} $ of the regular diffusion $Y$ is finite and $h(Y_{\tau_{l,r}})
= f(Y_{\tau_{l,r}})$, $\mathbb{P}^{\pi}$-almost surely for all $\pi
\in(0,1)$.
Then, by applying It\^{o}'s rule, we obtain
\begin{eqnarray*}
\mathbb{E}^{\pi} e^{- \lambda\tau_{l,r} } h(Y_{\tau_{l,r}}) &=&
\mathbb{E}^{\pi} e^{-
\lambda\tau_{l,r} } f(Y_{\tau_{l,r}}) \\
&=& f(\pi) + \mathbb{E}^{\pi} \int_0^{ \tau_{l,r} } e^{- \lambda t
} (-
\lambda+ \mathcal{A}_0 )f(Y_u)\, du \\
&=& f(\pi) - \mathbb{E}^{\pi} \int_0^{ \tau_{l,r} } e^{- \lambda t
} [
g(Y_u) + \lambda w(\mathbb{S}(Y_u)) ]\, du,
\end{eqnarray*}
and this shows $f(\cdot) = H_{l,r}[w](\cdot)$ on $(l,r)$.
\end{pf}

\begin{lemm}
\label{lemmexittimeexpectation2}
For $0< r <1 $, and $\tau_r := \inf\{ t \ge0 \dvtx Y_t \ge r \}$, let us define
\begin{eqnarray}
\label{defHr}
&&H_{r}[w](\pi) := \mathbb{E}^{\pi} \biggl[ \int_0^{\tau_{r}}
e^{-\lambda t }
[ g( Y_t ) + \lambda w( \mathbb{S}(Y_t) ) ] \,dt
\nonumber
\\[-8pt]
\\[-8pt]
\nonumber
&&\hspace*{156pt}{}+ e^{ -
\lambda\tau
_{r} } h( Y_{\tau_{r}} ) \biggr],\qquad  \pi\in[0,1],
\end{eqnarray}
which clearly equals $h(\cdot)$, for $\pi\ge r$.
For $\pi\in(0,r)$, the function $H_r[w](\cdot)$ can be computed by
taking the limit of \eqref{Hlrexplicit} as $l \searrow0$. That is,
\begin{eqnarray}
\label{Hrexplicit}
H_{r}[w](\pi) &=& \lim_{l \searrow0} H_{l,r}[w](\pi)\nonumber\\
&=&\psi(\pi) \biggl( - \frac{\eta(r)}{\psi(r)} \int_0^{r} u_2[w](y)
\,dy +
\int_{\pi}^{r} u_1[w](y) \,dy + \frac{h(r)}{\psi(r)}\biggr)\\
&&{}+ \eta(\pi) \int_0^{\pi} u_2[w](y)\, dy .\nonumber
\end{eqnarray}
The expression in \eqref{Hrexplicit} is twice-continuously
differentiable [on $(0,r)$] and solves
\begin{equation}
\label{differentialequationforH}
( - \lambda+ \mathcal{A}_0 )H_{r}[w](\pi) + g(\pi) + \lambda
w(\mathbb{S}(y)) =0.
\end{equation}
Moreover, the function $ H_r(\cdot) $ is continuous on $[0,1]$ with 
%
\begin{equation}
\label{continuityat0}
\lim_{\pi\searrow0} H_{r}[w](\pi) = w(p) = H_{r}[w](0).
\end{equation}
\end{lemm}

\begin{pf}
The point $\{0\}$ is a natural boundary for $Y$; therefore, we have
$\tau_r = \lim_{l \searrow0} \tau_{l,r}$, $\mathbb{P}^{\pi
}$-almost surely,
for $\pi\in(0,r)$. Then, the dominated convergence theorem (see
Remark \ref{remintegrabilitycond}) implies that $H_{r}[w](\pi) =
\lim
_{l \searrow0} H_{l,r}[w](\pi) $.

To compute the limit of $H_{l,r}[w](\pi)$ as $l \searrow0$, we first observe
\begin{eqnarray}
\label{limitofHterm}
&&\lim_{l \searrow0}
  h(l) \frac{\psi(\pi) \eta(r) - \psi(r) \eta(\pi) }{\psi
(l)\eta(r)-
\psi(r)\eta(l)}
+ h(r) \frac{\psi(l) \eta(\pi) - \psi(\pi) \eta(l) }{\psi
(l)\eta(r)-
\psi(r)\eta(l)}
\nonumber
\\[-8pt]
\\[-8pt]
\nonumber
&&\qquad= \psi(\pi) \frac{h(r)}{\psi(r)}.
\end{eqnarray}
Moreover, since $0 \le g(\cdot)+ \lambda w (\mathbb{S}(\cdot)) \le c
+ \lambda
\| w \|$,
we have
\begin{eqnarray*}
\int_{0}^{\pi} u_2[w](y)\, dy &\le&
\frac{c + \lambda\| w \|}{m_1 - m_2} \int_{0}^{\pi} \frac{2 \psi
(y)}{\sigma^2(y)} \,dy \\
&= &\frac{c + \lambda\| w \|}{m_1 - m_2} \int_{0}^{\pi}
\frac{\psi''(y)}{\lambda} \,dy\\
&=& \frac{(c/\lambda) + \| w \|}{m_1 - m_2} \psi'(\pi) < \infty\qquad \mbox{for $\pi< 1$,}
\end{eqnarray*}
and using \eqref{defuandC} we get
\begin{equation}
\label{limitsof_Cs}
\quad\lim_{l \searrow0} C_2 = \int_{0}^{r} u_2[w](y)\, dy
 \quad \mbox{and}\quad
\lim_{l \searrow0} C_1 =  \frac{\eta(r)}{\psi(r)} \int_{0}^{r}
u_2[w](y)\, dy .
\end{equation}
Finally letting $l \searrow0$ in \eqref{Hlrexplicit} and using the
limits found in \eqref{limitofHterm} and \eqref{limitsof_Cs}, we
obtain the expression in \eqref{Hrexplicit}. It is evident that this
expression 
is twice-continuously differentiable. Moreover, by direct computation
[using \eqref{defuandC}] it can be verified easily that it solves
the equation in \eqref{differentialequationforH}.

Clearly, $H_r[w](\cdot)$ is continuous on $(0,r)$ and $(r,1)$. The
continuity at $\{r\}$ can be checked by letting $\pi\nearrow r$ in the
expression given in
\eqref{Hrexplicit}, which goes to $h(r)$. 
To establish \eqref{continuityat0}, we first note that
$Y_t =0 $, for all $t > 0$, if $Y_0 = \pi=0$. 
This implies
\[
H_r[w](0) = \mathbb{E}^{\pi} \biggl[ \int_0^{\infty} e^{-\lambda t
} [ g(
0 ) + \lambda w( S(0) ) ] \,dt \biggr] = w(p).
\]
On the other hand, applying L'H\^{o}pital rule and using the explicit form
of $\psi(\cdot)$ and~$\eta(\cdot)$, we obtain
\begin{eqnarray*}
\lim_{\pi\searrow0} \psi(\pi) \int_{\pi}^{r} u_1[w] (y)\, dy&=&
\frac{2
\lambda  w(p) }{\mu^2 m_1 (m_1 - m_2)},
\\
\lim_{\pi\searrow0} \eta(\pi) \int_0^{\pi} u_2[w] (y)\, dy& =&
-\frac{2
\lambda  w(p) }{\mu^2 m_2 (m_1 - m_2)}.
\end{eqnarray*}
Since 
$m_1 \cdot m_2 = - 2 \mu^2 / \lambda$ [see \eqref
{quadraticequationform}], taking the limit 
in \eqref{Hrexplicit} gives $w(p)$, and this concludes the proof.
\end{pf}

Lemma~\ref{lemmexittimeexpectation2} shows that at the point $\pi=
r$, we have $(H_r[w])'(r+) =-1 $ and
\begin{eqnarray*}
(H_r[w])'(r-) &=&
\psi'(r) \biggl( - \frac{\eta(r)}{\psi(r)} \int_0^r u_2[w] (y) \,dy
+ \frac
{h(r)}{\psi(r)}\biggr)\\
&&{}+ \eta'(r) \int_0^r u_2[w] (y) \,dy .
\end{eqnarray*}
Since the Wronskian $\psi'(r) \eta(r) - \psi(r) \eta'(r)$ equals
$m_1 -
m_2$, we can rewrite the left derivative as
\[
(H_r[w])'(r-) =
\frac{1}{\psi(r)}
\biggl( - \int_0^r 2   \frac{ \psi(y) }{\sigma^2(y)}   [
g(y) +
\lambda w(\mathbb{S}(y))] \,dy
+ \psi'(r) h(r)
\biggr).
\]
Hence,
the derivative is continuous
at $\pi=r$ if and only if
\begin{eqnarray}\label{defB}
&&\hspace*{32pt}- \int_0^r 2   \frac{ \psi(y) }{\sigma^2(y)}   [ g(y) +
\lambda
w(\mathbb{S}(y))] \,dy
+ \psi'(r) h(r)
+ \psi(r) =0
\nonumber
\\[-8pt]
\\[-8pt]
\nonumber
&&\hspace*{32pt}\qquad \Longleftrightarrow\quad
B[w] (r) := \int_0^r \frac{ 2 \psi(y) }{\sigma^2(y)}   [ -
g(y) -
\lambda w(\mathbb{S}(y)) + \lambda h(y)]\, dy =0,
\end{eqnarray}
where the second equation follows after noting that
\[
\psi'(r) h(r)
+ \psi(r) = \int_0^r h (y) [ 2 \lambda\psi(y) / \sigma^2 (y)
]\,dy,
\]
which can be verified using $\lambda\psi(\cdot) = \mathcal{A}_0
\psi(\cdot)$.

\begin{lemm}
\label{lemmuniqueroot}
If $w(\cdot)$ is concave, then
the function $ - g(\pi) - \lambda w(\mathbb{S}(\pi)) + \lambda h(\pi
) $ has a
unique root $d[w] \in(0,1)$. The function $\pi\mapsto B[w] (\pi) $
equals zero for $\pi=0$, strictly increases on $(0,d[w])$ and strictly
decreases on $(d[w],1)$ with $\lim_{\pi\nearrow1} B[w] (\pi) = -
\infty$. Hence, there exists a unique point $r[w] \in(d[w],1)$ at
which $B[w] (r[w])=0$.
\end{lemm}

\begin{pf}
The function $\pi\mapsto-g(\pi) - \lambda w(\mathbb{S}(\pi)) +
\lambda h(\pi
) = c \pi- \lambda w(\mathbb{S}(\pi)) + \lambda(1-\pi) $ is convex and
continuous on $[0,1]$. At the point $\pi=0$,
it equals $\lambda( - w(p) + 1 ) \ge\lambda( - h(p) + 1 ) > 0 $;
and at
$\pi=1$, its value is $ -c - \lambda w(p) < 0$.
Hence, there exists a single point $d[w] \in(0,1)$ at which it is zero.
To the left of this point it is positive, and to the right it is negative.
Therefore, $\pi\mapsto B[w] (\pi)$ is zero at $\pi=0$, strictly
increases on $(0,d[w])$ and strictly decreases on $(d[w],1)$.
Also, observe that, for $\pi< 1$, 
\[
| B[w] (\pi) | \le(c + 2 \lambda) \int _0^{\pi}
\frac{ 2
\psi(y) }{\sigma^2(y)}\,dy = \frac{c + 2 \lambda}{\lambda}
\int _0^{\pi} \psi''(y) \,dy = \frac{c + 2 \lambda}{\lambda}
\psi'(\pi) < \infty
\]
and
\begin{eqnarray*}
&&\int_{d[w] + \delta}^{1}
\frac{2 \psi(y) }{\sigma^2(y)}   [ - g(y) - \lambda w(\mathbb
{S}(y)) +
\lambda h(y)] \,dy
\\
&&\qquad\le\Bigl(
\min_{\pi\in[d[w] + \delta, 1] } \{ - g(\pi) - \lambda w(\mathbb
{S}(\pi)) +
\lambda h(\pi) \}
\Bigr)
\cdot
\int_{d[w] + \delta}^{1} \frac{2 \psi(y) }{\sigma^2(y)}\, dy \\
&&\qquad= \Bigl(
\min_{\pi\in[d[w] + \delta, 1] } \{ - g(\pi) - \lambda w(\mathbb
{S}(\pi)) +
\lambda h(\pi) \}
\Bigr)
\cdot[\psi'(\pi)]|^{1}_{d[w] + \delta} = - \infty
\end{eqnarray*}
for all $\delta\in(0, 1- d[w])$, where the last equality follows
using \eqref{defpsiandeta}. Hence, we conclude that $B[w](\pi) $
goes to $-\infty$ as $\pi\to1$, and this implies that it has a unique
root $r[w] \in( d[w] ,1 )$.
\end{pf}

\begin{rem}
\label{remboundsonr}
For two concave functions $w_1(\cdot) \le w_2 (\cdot)$ satisfying
Assumption \ref{assa1}, we have $ B[w_1](\cdot) \ge B[w_2](\cdot) $; therefore
$r[w_1] \ge r[w_2]$.
If we select the zero function (which equals zero on $[0,1]$),
direct computation yields
\[
B[0](\pi) = \frac{\psi(\pi)}{\pi(1-\pi)} \biggl[ - \pi\biggl((m_1-1)
\frac{c}{\lambda} +m_1 \biggr) +m_1 \biggr],
\]
and for $h(\pi)=1-\pi$, we get
\[
B[h](\pi) = \frac{\psi(\pi)}{\pi(1-\pi)} \biggl[ - \pi\biggl((m_1-1)
\frac{c}{\lambda} +m_1 p \biggr) + m_1 p \biggr].
\]
Hence,
we have the bounds
\begin{equation}
\label{boundsonr}
\frac{m_1 p}{ (m_1-1) ({c}/{\lambda}) + m_1 p}
\le
r[w]
\le
\frac{m_1}{ (m_1-1) ({c}/{\lambda}) + m_1}.
\end{equation}
\end{rem}

Observe that, if the function $w(\cdot)$ is concave,
$H_{r[w]}[w](\cdot)$ is
continuously differentiable on $(0,1)$. On $(r[w],1)$,
$H_{r[w]}[w](\cdot) $ coincides with $ h(\cdot)$, and 
\[
(- \lambda+ \mathcal{A}_0 ) H_{r[w]}[w](\pi) + g(\pi) + \lambda
w(\mathbb{S}(\pi)) =
- \lambda h(\pi) + g(\pi) + \lambda w(\mathbb{S}(\pi)) > 0,
\]
since $d[w] < r[w]$.

On $(0,r[w])$, the function $H_{r[w]}[w](\cdot)$ solves
$(- \lambda+ \mathcal{A}_0 ) H_{r[w]}[w](\pi) + g(\pi) + \lambda
w(\mathbb{S}(\pi))
=0 $.
In Appendix \hyperref[secappendixB]{B}, we also show that
\begin{equation}
\label{inequalityonstrictconcavity}
\lambda H_{r[w]}[w](\pi) - g(\pi) - \lambda w(\mathbb{S}(\pi)) < 0\qquad
 \mbox
{for $0 < \pi< r[w] $.}
\end{equation}
Since $\mathcal{A}_0 H_{r[w]}[w](\pi) = ( \sigma^2(\pi) / 2) \cdot(
H_{r[w]}[w](\pi))''$,
the inequality in \eqref{inequalityonstrictconcavity} implies that
$H_{r[w]}[w](\cdot)$ is strictly concave and $H_{r[w]}[w](\cdot) <
h(\cdot)$ on $(0,r[w])$.

Finally, the (strict) concavity on $(0,r[w])$ and the ``smooth-fit'' at
$\{r[w]\}$ imply that $H_{r[w]}[w](\cdot)$ is also concave on $(0,1)$.
The following remark is a summary of the analytical properties of
$H_{r[w]}[w](\cdot)$ described above.

\begin{rem}
\label{remanalyticalpropertiesofHr}
Suppose that the function $w(\cdot)$ is concave. Then,
$H_{r[w]}[w](\cdot)$ is nonnegative, continuous and concave on $[0,1]$.
It is continuously differentiable on $(0,1)$, twice-continuously
differentiable on $(0,1) \setminus\{r[w]\}$, and it satisfies
\begin{eqnarray}
\label{variationalinequalititesforHr}
\qquad&&\left\{\matrix{
H_{r[w]}[w](\pi)= h(\pi)\vspace*{2pt}\cr
(- \lambda+ \mathcal{A}_0 ) H_{r[w]}[w](\pi) + g(\pi) + \lambda
w(\mathbb{S}(\pi)) >0}\right\},\qquad
 \pi\in(r[w],1),\hspace*{-14pt}
\nonumber
\\[-8pt]
\\[-8pt]
\nonumber
\qquad &&\left\{\matrix{
H_{r[w]}[w](\pi)< h(\pi)\vspace*{2pt}\cr
(- \lambda+ \mathcal{A}_0 ) H_{r[w]}[w](\pi) + g(\pi) + \lambda
w(\mathbb{S}(\pi)) =0}\right\},\qquad
 \pi\in(0,r[w]).\hspace*{-14pt}
\end{eqnarray}

\end{rem}

\begin{lemm}
\label{lemmJequalsHr}
If $w(\cdot)$ is concave, we have $J[w](\cdot) = H_{r[w]}[w](\cdot)$,
and $\tau_{r[w]} := \inf\{t \ge0\dvtx Y_t \ge r[w] \}$ is an optimal
stopping time for \eqref{Jreloaded}. 
\end{lemm}
\begin{pf}
For $\pi\in(0,1)$, let $\tau$ be an $\mathbb{F}^Y$-stopping time,
and $\tau
_{l,r}$ be the exit time of $Y$ from $(l,r)$ for $0< l \le r < 1$.
Then, by It\^{o}'s rule
\begin{eqnarray*}
&&e^{- \lambda\cdot\tau\wedge\tau_{l,r} } H_{r[w]}[w]( Y_{\tau
\wedge\tau_{l,r} } )\\
&&\qquad= H_{r[w]}[w](\pi)
+ \int_0^{\tau\wedge\tau_{l,r}} e^{- \lambda t} (-\lambda+
\mathcal{A}
_0)H_{r[w]}[w](Y_t)\, dt\,\\
&&\qquad\quad {}+ \int_0^{\tau\wedge\tau_{l,r}} e^{- \lambda t} \sigma(Y_t) \bigl(
H_{r[w]}[w] \bigr)' (Y_t)\, d\widehat{W}_t.
\end{eqnarray*}
The function $H_{r[w]}[w](\cdot)$ is continuously differentiable on
$(0,1)$. Its derivative is therefore bounded on $[l,r]$ and
$\| \sigma(\cdot) \| \le| \mu| $. Then, taking expectations above gives
\begin{eqnarray*}
&&\mathbb{E}^{\pi} e^{- \lambda\cdot\tau\wedge\tau_{l,r} }
H_{r[w]}[w]( Y_{\tau\wedge\tau_{l,r} } ) \\
&&\qquad= H_{r[w]}[w](\pi)
+\mathbb{E}^{\pi} \int_0^{\tau\wedge\tau_{l,r}} e^{- \lambda t}
(-\lambda+
\mathcal{A}_0)H_{r[w]}[w](Y_t)\, dt \\
&&\qquad\ge H_{r[w]}[w](\pi) -\mathbb{E}^{\pi} \int_0^{\tau\wedge\tau
_{l,r}} e^{-
\lambda t} \bigl( g(Y_t) + \lambda w(\mathbb{S}(Y_t)) \bigr)\, dt,
\end{eqnarray*}
where the inequality is due to
\eqref{variationalinequalititesforHr}. Since both boundaries are
natural, we first let $r \nearrow1$ and then $l \searrow0$ to
obtain
\begin{eqnarray}
\qquad \mathbb{E}^{\pi} e^{- \lambda\cdot\tau} h ( Y_{\tau}
) &\ge&
\mathbb{E}^{\pi} e^{- \lambda\cdot\tau} H_{r[w]}[w]( Y_{\tau}
)
\nonumber
\\[-8pt]
\\[-8pt]
\nonumber
&\ge&H_{r[w]}[w](\pi) -\mathbb{E}^{\pi} \int_0^{\tau} e^{-
\lambda t} \bigl( g(Y_t) + \lambda w(\mathbb{S}(Y_t)) \bigr)\, dt
\end{eqnarray}
thanks to dominated convergence theorem (see
Remark~\ref{remintegrabilitycond}), and this shows $
H_{r[w]}[w](\cdot) \le J[w](\cdot)$ on $(0,1)$.

When we repeat the steps above with
$\tau= \tau_{r[w]} = \inf\{ t \ge0\dvtx Y_t \ge r[w]\}$,
the inequalities become equalities again by \eqref{variationalinequalititesforHr}. Hence, $J[w](\cdot) =
H_{r[w]}[w](\cdot)$ on $(0,1)$.

If $Y_0 = \pi= 0$; then $Y_t = 0$, for $t \ge0$, and
\begin{eqnarray*}
J[w] (0) &= &
\inf_{ \tau}
\mathbb{E}^{\pi} [ \lambda w(\mathbb{S}(0)) ( 1- e^{- \lambda
\tau}) + e^{-
\lambda\tau} h (0) ]\\
& = & w(\mathbb{S}(0)) = w(p) = H_{r[w]}(0),
\end{eqnarray*}
thanks to Lemma~\ref{lemmexittimeexpectation2} [note that $w(p)
\le
h(p) < h(0)$]. Moreover, this value is attained by selecting $\tau=
\infty= \tau_{r[w]}$. Similarly, if $Y_0 =1 $, we have $ J[w] (1) =
h(1)=0 = H_{r[w]}(1)$, which is attained by $\tau= 0= \tau_{r[w]}$.
\end{pf}

\section{The value function and an optimal detection rule}
\label{secsolution}
Using the dynamic programming operator $J$, let us define the sequence
of functions
\[
v_0 \equiv h(\cdot) \quad   \mbox{and} \quad   v_{n+1} (\cdot) :=
J[v_n] (\cdot) = H_{r[v_n]}[v_n](\cdot)\qquad
\mbox{for $n \in\mathbb{N}$}.
\]

\begin{rem}
\label{vnsconverge}
The sequence $(v_n )_{n \in\mathbb{N}}$ is nonincreasing, and each
element of
the sequence is a nonnegative, continuous and concave function on
$[0,1]$. 
\end{rem}
\begin{pf}
We have $v_1 (\cdot) = J[h] (\cdot) \le h (\cdot)= v_0(\cdot)$,
where the inequality follows from the definition of the operator $J$ in
\eqref{Jreloaded}. Next, assume that $v_{n}(\cdot) \le v_{n-1}
(\cdot
)$, for some $n \in\mathbb{N}$. Then Remark~\ref
{rempropertiesofJ} implies
$v_{n+1}(\cdot)= J[v_{n}](\cdot) \le\break J[v_{n-1}] (\cdot) =
v_{n}(\cdot)$, and this shows that the sequence $(v_n )_{n \in\mathbb
{N}}$ is
nonincreasing by induction. Finally, since $v_0 (\cdot) =h(\cdot)$ is
nonnegative, continuous and concave,
these properties also hold for each $v_n (\cdot)$, $n \in\mathbb
{N}$, by
induction thanks to Remark~\ref{remanalyticalpropertiesofHr}.
\end{pf}

For $n \in\mathbb{N}$, let $\pi_n := r[v_{n-1}]$ be the solution of the
equation $B[v_{n-1}](r)=0$ [see~\eqref{defB}]. Since $v_{n-1}(\cdot)$
is concave and satisfies Assumption \ref{assa1}, this equation has a unique root
on $(0,1)$ thanks to Lemma~\ref{lemmuniqueroot}. Moreover,
Remark~\ref
{remanalyticalpropertiesofHr} and Lemma~\ref{lemmJequalsHr}
imply that $v_n(\cdot)$ is continuously differentiable on $(0,1)$,
twice-continuously differentiable on $(0,1) \setminus\{ \pi_n \}$ and
solves the variational inequalities\looseness=1
\begin{eqnarray*}
&&\left\{\matrix{
v_n(\pi)= h(\pi) \cr
(- \lambda+ \mathcal{A}_0 ) v_n (\pi) + g(\pi) + \lambda
v_{n-1}(\mathbb{S}(\pi)) &>0}\right\},\qquad
 \pi\in(\pi_n,1),
\\
&&\left\{\matrix{
v_n(\pi)< h(\pi) \cr
(- \lambda+ \mathcal{A}_0 ) v_n(\pi) + g(\pi) + \lambda
v_{n-1}(\mathbb{S}(\pi)) =0}\right\},\qquad
 \pi\in(0,\pi_n).
\end{eqnarray*}
Observe that $\pi_n = \inf\{ \pi\in[0,1] \dvtx v_n (\pi) = h(\pi) \}
$; hence,
$\{ \pi_n \}_{ n \in\mathbb{N}}$ is nondecreasing.

Let $v_{\infty}(\cdot) := \inf_{n \in\mathbb{N}} v_{n}(\cdot)$ be the
pointwise limit of 
$( v_n )_{n \in\mathbb{N}}$. Then
dominated convergence theorem gives
\begin{eqnarray}
\label{fixedpointofJ1}
v_{\infty}(\pi) &=& \inf_{n \in\mathbb{N}} J[v_{n-1}] (\pi) \nonumber\\
&=& \inf_{n \in\mathbb{N}} \inf_{ \tau\in\mathbb{F}^{Y} }
\mathbb{E}^{\pi} \biggl[ \int_0^{\tau} e^{-\lambda t } [ g(
Y_t ) + \lambda
v_{n-1}( \mathbb{S}(Y_t) ) ] \,dt + e^{ - \lambda\tau} h( \Pi
^{\pi}_{\tau
} ) \biggr]
\nonumber\\
&=& \inf_{ \tau\in\mathbb{F}^{Y} } \inf_{n \in\mathbb{N}}
\mathbb{E}^{\pi} \biggl[ \int_0^{\tau} e^{-\lambda t } [ g(
Y_t ) + \lambda
v_{n-1}( \mathbb{S}(Y_t) ) ]\, dt + e^{ - \lambda\tau} h( \Pi
^{\pi}_{\tau
} ) \biggr] \\
&= &\inf_{ \tau\in\mathbb{F}^{Y} }
\mathbb{E}^{\pi} \biggl[ \int_0^{\tau} e^{-\lambda t } [ g(
Y_t ) + \lambda
v_{\infty}( \mathbb{S}(Y_t) ) ]\, dt + e^{ - \lambda\tau} h(
\Pi^{\pi
}_{\tau} ) \biggr]\nonumber\\
&=&J[v_{\infty}] (\pi),\nonumber
\end{eqnarray}
which shows that the function $v_{\infty}(\cdot)$ is a fixed point of
the operator $J$.

\begin{lemm}
\label{lemmuniformconvergence}
The sequence $(v_n)_{n \ge1}$ converges to $v_{\infty}(\cdot)$
uniformly on $[0,1]$. More precisely, we have
\begin{equation}
\label{equniformconvergence}
v_{\infty}(\pi) \le v_{n}(\pi) \le v_{\infty}(\pi) + (1-p)^n
(1-\pi)
\end{equation}
for all $n \in\mathbb{N}$.
\end{lemm}
\begin{pf}
The first inequality in \eqref{equniformconvergence} is immediate
since the sequence $(v_n)_{n \in\mathbb{N}}$ is nonincreasing. The second
inequality is also obvious for $n=0$ as $v_{\infty} (\cdot) = \inf_{n
\in\mathbb{N}} v_n (\cdot) \ge0$. Assume the second inequality
holds for some
$n \in\mathbb{N}$. This implies that $v_n(\mathbb{S}(\pi)) = v_n
(\pi+p (1-\pi)) \le
v_{\infty} (\pi+p (1-\pi)) + (1-p)^n (1 - \pi- p (1-\pi)) =
v_{\infty
}(\mathbb{S}(\pi)) + (1-p)^{n+1} (1-\pi)$. Then we have\vspace*{1pt}
\begin{eqnarray*}
v_{n+1} (\pi) &=& J[v_n] (\pi)\\[1pt]
&=& \inf_{ \tau\in\mathbb{F}^{Y} }
\mathbb{E}^{\pi} \biggl[ \int_0^{\tau} e^{-\lambda t } [ g(
Y_t ) + \lambda
v_{n}( \mathbb{S}(Y_t) ) ]\, dt + e^{ - \lambda\tau} h( \Pi
^{\pi}_{\tau}
) \biggr]\\[1pt]
&\le&\inf_{ \tau\in\mathbb{F}^{Y} }
\mathbb{E}^{\pi} \biggl[ \int_0^{\tau} e^{-\lambda t } [ g(
Y_t ) + \lambda
v_{\infty}( \mathbb{S}(Y_t) ) + \lambda(1-p)^{n+1} (1-Y_t )]\, dt\\[1pt]
&&\hspace*{226pt}{} + e^{ -
\lambda\tau} h( \Pi^{\pi}_{\tau} ) \biggr]
\\[1pt]
&\le&
\inf_{ \tau\in\mathbb{F}^{Y} }
\mathbb{E}^{\pi} \biggl[ \int_0^{\tau} e^{-\lambda t } [ g(
Y_t ) + \lambda
v_{\infty}( \mathbb{S}(Y_t) ) ] \,dt + e^{ - \lambda\tau} h(
\Pi^{\pi
}_{\tau} ) \biggr]
\\[1pt]
& &{}       + \mathbb
{E}^{\pi} \biggl[ \int
_0^{\infty} e^{-\lambda t } \lambda(1-p)^{n+1} (1-Y_t )\, dt \biggr].
\end{eqnarray*}
Since $v_{\infty}(\cdot)$ satisfies $v_{\infty}(\cdot) =
J[v_{\infty
}](\cdot)$, the last inequality gives
\begin{eqnarray*}
v_{n+1} (\pi) &\le& v_{\infty} (\pi) +
\int_0^{\infty} e^{-\lambda t } \lambda(1-p)^{n+1} \mathbb{E}^{\pi
}
[1-Y_t ] \,dt
\\
&=& v_{\infty} (\pi) + (1-p)^{n+1} (1-\pi),
\end{eqnarray*}
where we used the martingale property of $Y$ to justify the last
equality. This\break shows
the second inequality in \eqref{equniformconvergence} for $n+1$, and
the proof is complete by induction. 
\end{pf}

\begin{cor}
\label{corvinftysatisfytheassumptions}
The uniform convergence in Lemma~\textup{\ref{lemmuniformconvergence}}
implies that $v_{\infty}(\cdot)$ is continuous on $[0,1]$. Moreover, as
the infimum of nonnegative concave functions $v_n(\cdot)$\textup{'}s, it is
also nonnegative and concave.
\end{cor}

Corollary~\ref{corvinftysatisfytheassumptions} and the identity
$ v_{\infty}(\cdot) = J[v_{\infty}](\cdot)$ [see (\ref
{fixedpointofJ1})] allow us to conclude that $v_{\infty}(\cdot)$ is
continuously differentiable on $(0,1)$ and twice-continuously
differentiable on $(0,1) \setminus\{ \pi_{\infty} \}$, where $\pi
_{\infty} := r [v_{\infty}]$ is the unique root 
of the equation $B[v_{\infty}](r) =0$ defined in \eqref{defB}.
Furthermore, $v_{\infty}(\cdot)$ satisfies
\begin{eqnarray}
\label{variationalinequalititesforvinfty}
&&\left\{\matrix{
v_{\infty}(\pi)= h(\pi) \cr
(- \lambda+ \mathcal{A}_0 ) v_{\infty} (\pi) + g(\pi) + \lambda
v_{{\infty
}}(\mathbb{S}(\pi)) >0}\right\},\qquad
 \pi\in(\pi_{\infty},1),
\nonumber
\\[-8pt]
\\[-8pt]
\nonumber
&&\left\{\matrix{
v_{\infty}(\pi)< h(\pi) \cr
(- \lambda+ \mathcal{A}_0 ) v_{\infty}(\pi) + g(\pi) + \lambda
v_{\infty}(\mathbb{S}
(\pi)) =0}\right\},\qquad
 \pi\in(0 ,\pi_{\infty}),
\end{eqnarray}
which also implies that $\pi_{\infty} = \inf\{ \pi\in[0,1] \dvtx
v_{\infty} (\pi) = h(\pi) \}$. Since $v_n (\cdot) \searrow v(\cdot)$,
we have $ \pi_n \nearrow\pi_{\infty}$.

\begin{prop}
\label{propVandoptimaldetectiontime}
The function $v_{\infty}(\cdot)$ is the value function $V(\cdot)$ of
the optimal stopping problem in \eqref{defV}, and the first entrance
time $\widetilde{\tau}_{\pi_{\infty}}$ of the process $\Pi$ to the
interval $[\pi_{\infty},1] $ is an optimal solution for the
change-detection problem in~\eqref{defU}.
\end{prop}

\begin{pf}
The claim is obvious if $\Pi_0 = \pi= 1$; both $v_{\infty}(1)$ and
$V(1)$ are nonnegative and bounded by $h(1)=0$, which is also the
expected reward in \eqref{defV} by stopping immediately.

For $\pi\in(0,1)$ 
and $0 < l \le r < 1$, let $\widetilde{\tau}_{[0,l]}$ and $\widetilde
{\tau}_{[r,1]}$ be respectively, the entrance times of the process
$\Pi
$ to the intervals $[0,l]$ and $[r,1]$.
Also, define $\widetilde{\tau}_{l,r} : = \widetilde{\tau}_{[0,l]}
\wedge\widetilde{\tau}_{[r,1]}$. Then for an $\mathbb{F}$-stopping
time $\tau$,
It\^{o}'s rule gives
\begin{eqnarray*}
v_{\infty} (\Pi_{ \tau\wedge\widetilde{\tau}_{l,r} })
&=& v_{\infty} (\pi) +
\int_0^{ \tau\wedge\widetilde{\tau}_{l,r} } [ (- \lambda+
\mathcal{A}
_0) v_{\infty} (\Pi_{u-})
+ \lambda v_{\infty} (\mathbb{S}(\Pi_{u-})) ] \,du\\
&&{}+ \int_0^{ \tau\wedge\widetilde{\tau}_{l,r} } \mu\Pi_{u-} (1-
\Pi
_{u-}) v'_{\infty} (\Pi_{u-})\, d \widehat{W}_u \\
&&{}+ \int_0^{ \tau\wedge\widetilde{\tau}_{l,r} } [ v_{\infty}
(\mathbb{S}(\Pi_{u-})) -
v_{\infty} (\Pi_{u-}) ] ( d N_u - \lambda \,du ).
\end{eqnarray*}
Since the function $v_{\infty}$ is bounded, the stochastic integral
with respect to the martingale $\{ N_t - \lambda t \}_{t \ge0}$ is a
square-integrable martingale stopped at $\tau\wedge\widetilde{\tau
}_{l,r}$ [whose expectation is finite due to \eqref
{exittimeexpectationuniformlybounded}]. Similarly, so is the
integral with respect to~$\widehat{W}$ as $v'_{\infty}$ is continuous
and bounded on $[l,r]$. Then taking expectations, we obtain\looseness=1
\begin{eqnarray}
\label{takingexpectations1}
&&\mathbb{E}^{\pi} v_{\infty} (\Pi_{ \tau
\wedge\widetilde{\tau}_{l,r} })\nonumber\\
&&\qquad= v_{\infty} (\pi) + \mathbb{E}^{\pi} \int_0^{ \tau\wedge
\widetilde{\tau}_{l,r} } [ (- \lambda+ \mathcal{A}_0)
v_{\infty}
(\Pi_{u-})
+ \lambda v_{\infty} (\mathbb{S}(\Pi_{u-})) ]\, du
\nonumber
\\[-8pt]
\\[-8pt]
\nonumber
&&\qquad\ge v_{\infty} (\pi) - \mathbb{E}^{\pi} \int_0^{ \tau\wedge
\widetilde{\tau}_{l,r} } g(\Pi_{u-} )\, du\\
&&\qquad = v_{\infty} (\pi) -
\mathbb{E}^{\pi} \int_0^{ \tau\wedge\widetilde{\tau}_{l,r} }
g(\Pi_{u}
) \,du\nonumber
\end{eqnarray}
thanks to the inequalities in \eqref{variationalinequalititesforvinfty}.

The left boundary $\{0\}$ is natural for the diffusion in \eqref
{dynamicsofY}. Between two arrivals of $N$, the process $\Pi$ follows
these dynamics, and at an arrival time $T_n$ it jumps to the right by\
an amount of $p (1- \Pi_{T_n} )$. Hence, as we let $l \searrow0$,
$\widetilde{\tau}_{[0,l]}$ goes to~$\infty$ thanks to strong Markov
property, and 
$ \widetilde{\tau}_{l,r} \nearrow\widetilde{\tau}_{[r,1]}$. Moreover,
$\lim_{t \to\infty} \Pi_t = 1$, and $\Pi_t < 1$ [since $\Phi_t <
\infty
$ in \eqref{eqPiandPhi}] for finite $t$, if $\pi< 1$. Hence, as $r
\nearrow1$, we have $\widetilde{\tau}_{[r,1]} \nearrow\infty$.
Therefore, when we let $l \searrow0$ and $r \nearrow1$ in \eqref
{takingexpectations1},
bounded convergence and monotone convergence theorems give
\begin{equation}
\label{expectationsinequality}
\mathbb{E}^{\pi} v_{\infty} (\Pi_{ \tau})
\ge
v_{\infty} (\pi) - \mathbb{E}^{\pi} \int_0^{ \tau} g(\Pi_{u} )\,du .
\end{equation}
Also note that we have $\mathbb{E}^{\pi} h (\Pi_{ \tau}) \ge
\mathbb{E}^{\pi}
v_{\infty} (\Pi_{ \tau})$. Then we obtain
\[
\mathbb{E}^{\pi} \biggl[ \int_0^{ \tau} g(\Pi_{u} )\, du +h (\Pi_{
\tau})
\biggr]
\ge
v_{\infty} (\pi) ,
\]
which implies that $v_{\infty} (\pi) \le V (\pi) $ on $(0,1)$.

When we replace $\tau$ in \eqref{takingexpectations1} with the
entrance time $\widetilde{\tau}_{\pi_{\infty}}$, the inequality in
\eqref{expectationsinequality} becomes an equality. Then the equality
$\mathbb{E}^{\pi} h (\Pi_{ \widetilde{\tau}_{\pi_{\infty}} }) =
\mathbb{E}^{\pi}
v_{\infty} (\Pi_{ \widetilde{\tau}_{\pi_{\infty}} })$ yields
\begin{equation}
\label{Vattained}
v_{\infty}(\pi) = \mathbb{E}^{\pi} \biggl[\int_0^{ \widetilde
{\tau}_{[\pi
_{\infty},1 ]} } g(\Pi_{u} )\, du
+ h \bigl(\Pi_{ \widetilde{\tau}_{[\pi_{\infty},1 ]} }\bigr) \biggr],
\end{equation}
and this implies $V(\pi)= v_{\infty}(\pi)$, for $\pi\in(0,1)$.

To show the same equality for $\Pi_0 = \pi= 0 $, we first note that
$\Pi_t = 0 $ for $t < T_1$, and $\Pi_{T_1} = p $ if the process $\Pi$
starts from the point $\{0\}$.
Also note that the identity $v_{\infty}(\cdot) = J[v_{\infty}](\cdot
)$ implies
$v_{\infty}(0)= v_{\infty}(p)$ [see \eqref{defJ}].
Then for an $\mathbb{F}$-stopping time $\tau$, by modifying the arguments
above, we get
\begin{eqnarray}
\label{takingexpectations2}
&&\hspace*{-4pt}\mathbb{E}^{0} v_{\infty} \bigl(\Pi_{ \tau\wedge( \widetilde{\tau
}_{l,r} \circ   \theta_{T_1} ) }\bigr)\nonumber\\
&&\hspace*{-4pt}\qquad = v_{\infty}(0) + \mathbb{E}^{0} 1_{ \{ \tau\ge T_1 \} } \int_{T_1}^{ \tau\wedge(
\widetilde{\tau}_{l,r} \circ   \theta_{T_1} ) } [ (-
\lambda+ \mathcal{A}_0) v_{\infty} (\Pi_{u})\nonumber\\
&&\hspace*{181pt}{} + \lambda v_{\infty}
(\mathbb{S}(\Pi_{u-})) ]\, du
\\
 &&\hspace*{-4pt}\qquad\ge
v_{\infty} (0)
- \mathbb{E}^{0} 1_{ \{ \tau\ge T_1 \} } \int_{T_1}^{ \tau\wedge(
\widetilde
{\tau}_{l,r} \circ   \theta_{T_1} ) } g (\Pi_{u})\, du\nonumber\\
&&\hspace*{-4pt}\qquad=
v_{\infty} (0)
- \mathbb{E}^{0} \int_{0}^{ \tau\wedge( \widetilde{\tau}_{l,r}
\circ
\theta_{T_1} ) } g (\Pi_{u}) \,du,\nonumber
\end{eqnarray}
where $\theta$ is the time-shift operator.
Letting $l \searrow0$ and $r \nearrow1$ in \eqref{takingexpectations2},
and using the inequality
$\mathbb{E}^{0} h (\Pi_{ \tau})
\ge\mathbb{E}^{0} v_{\infty} (\Pi_{ \tau})
$,
we obtain $v_{\infty}(0) \le V(0)$.

Replacing $\tau$ above
with
$\widetilde{\tau}_{[\pi_{\infty},1 ]}$, we get
equalities in \eqref{takingexpectations2}.\vspace*{1pt} Then, letting $l \searrow
0$, $r \nearrow1$, and using
the equality
$\mathbb{E}^{0} h (\Pi_{ \widetilde{\tau}_{[\pi_{\infty},1 ]} })
= \mathbb{E}^{0} v_{\infty} (\Pi_{ \widetilde{\tau}_{[\pi_{\infty
},1 ]} })
$ we obtain \eqref{Vattained} for $\pi=0$. Hence, we have
$v_{\infty}(0) = V(0)$, and this concludes the proof.
\end{pf}

%

\begin{rem}
For $ \varepsilon> 0$, let us fix $n \in\mathbb{N}$ such that $n \ge
\ln
(\varepsilon) / \ln(1-p)$
and $ v_{n} (\cdot) \le v_{\infty} (\cdot) + \varepsilon$, on $[0,1]$.
The exit time $\widetilde{\tau}_{\pi_n}$ of $\Pi$
from the interval $[0,\pi_n) = \{ \pi\in[0,1] \dvtx v_n(\pi) < h(\pi)\}$
is $\varepsilon$-optimal
for the problem in \eqref{defV}. That is,
\begin{equation}
\label{ineqforepsilonoptimaltime}
\mathbb{E}^{\pi} \biggl[ \int_{0}^{\widetilde{\tau}_{\pi_n}}
g(\Pi_t)\, dt
+ h( \Pi_{\widetilde{\tau}_{\pi_n}} ) \biggr] \le V
(\pi) +
\varepsilon\qquad
 \mbox{for all $\pi\in[0,1]$.}
\end{equation}
\end{rem}
\begin{pf}
For $\pi> 0$, 
a localization argument and It\^{o}'s rule (as in the proof of
Proposition~\ref{propVandoptimaldetectiontime}) give
\begin{eqnarray*}
\mathbb{E}^{\pi} v_{\infty} (\Pi_{\widetilde{\tau}_{\pi_n}} )
&=& v_{\infty}(\pi) + \mathbb{E}^{\pi} \int_0^{ \widetilde{\tau
}_{\pi_{n}} }
[ (- \lambda+ \mathcal{A}_0) v_{\infty} (\Pi_{u-})
+ \lambda v_{\infty} (\mathbb{S}(\Pi_{u-})) ] \,du \\
&=&
v_{\infty}(\pi) - \mathbb{E}^{\pi} \int_0^{ \widetilde{\tau
}_{\pi_{n},1 } }
g(\Pi_{u} ) \,du ,
\end{eqnarray*}
where the last equality follows from \eqref{variationalinequalititesforvinfty} (recall that $\pi_n \le\pi
_{\infty}$). Note that $\widetilde{\tau}_{ \pi_{n} } < \infty$ and
$v_{n} (\Pi_{\widetilde{\tau}_{\pi_n}} ) = h (\Pi_{\widetilde
{\tau}_{\pi
_n}} )$, $\mathbb{P}^{\pi}$-almost surely. Then the inequality $
v_{n} (\cdot)
\le v_{\infty} (\cdot) + \varepsilon$ yields
\begin{eqnarray*}
\mathbb{E}^{\pi} h (\Pi_{\widetilde{\tau}_{\pi_n}} ) -
\varepsilon&=&
\mathbb{E}^{\pi} v_{n} (\Pi_{\widetilde{\tau}_{\pi_n}} ) -
\varepsilon
\\
&\le&\mathbb{E}^{\pi} v_{\infty} (\Pi_{\widetilde{\tau}_{\pi_n}} )\\
&=&
v_{\infty}(\pi) - \mathbb{E}^{\pi} \int_0^{ \widetilde{\tau
}_{\pi_{n}} }
g(\Pi_{u} ) \,du ,
\end{eqnarray*}
and \eqref{ineqforepsilonoptimaltime} follows.

For $\pi= 0$, we have
\begin{eqnarray*}
\mathbb{E}^{0} \int_0^{ \widetilde{\tau}_{\pi_{n}} }
g(\Pi_{u} )\, du + h (\Pi_{\widetilde{\tau}_{\pi_n}} )
&=& \mathbb{E}^{0} \int_{T_1}^{ \widetilde{\tau}_{\pi_{n}} \circ
  \theta
_{T_1} }
g(\Pi_{u} )\, du + h (\Pi_{\widetilde{\tau}_{\pi_n} \circ
\theta
_{T_1}} )
\\
&\le& V(p) + \varepsilon= V(0) + \varepsilon,
\end{eqnarray*}
where the inequality is due the strong Markov property (and also the
result already proved above for $\pi=p > 0$), and
the last equality follows from the identity $V(0) = J[V](0)= V(p)$.
\end{pf}
\section{Variational formulation}
\label{secvariational}

In this section, we solve the variational formulation of the problem
where the objective is to minimize the expected detection delay
$\mathbb{E}^{\pi
} (\tau- \Theta)^+ $ over all $\mathbb{F}$-stopping times for which
the false
alarm probability $ \mathbb{P}^{\pi} (\tau< \Theta) $ is less than
or equal
to some predetermined value $\alpha\in(0,1)$. The optimality of $\tau
=0$ is immediate when $\pi=1$; hence, this case is excluded below.

When $\pi\in(0,1)$, $\tau= 0$ is also an optimal solution if $\alpha
\ge1- \pi$. On the other hand, if $\pi=0$ and $\alpha\ge1-p $, the
first arrival time $T_1$ of $N$ yields a false alarm probability of
$1-p$ and its expected delay is still zero [see (\ref
{distributionofzeta})--(\ref{defTheta})]. 

If none of these trivial cases hold, we can find an optimal stopping
time (for the variational formulation) using the solution of the
problem in \eqref{defU} as explained in \cite{Sh78}. More precisely,
let $\pi_{\infty}(c)$ be the optimal threshold found in Section~\ref
{secsolution} as a function of $c$, and let $\widetilde{\tau}_{\pi
_{\infty}(c)}$ be the corresponding exit time of the process~$\Pi$. For
a given value of $\alpha$,
assume there exists a value of $c > 0 $ such that the false alarm
probability $\mathbb{P}^{\pi} ( \widetilde{\tau}_{\pi_{\infty}(c)
} < \Theta)
= \mathbb{E}^{\pi} h ( \Pi_{ \widetilde{\tau}_{\pi_{\infty
}(c) } } )$
equals $ \alpha$. Then $\widetilde{\tau}_{\pi_{\infty}(c)}$ solves the
variational formulation. Indeed, the optimality of $\widetilde{\tau
}_{\pi_{\infty}(c)}$ for the original problem in \eqref{defU} implies
that, for any $\mathbb{F}$-stopping time $\tau$, we have
\[
c   \mathbb{E}^{\pi} \bigl(\widetilde{\tau}_{\pi_{\infty}(c)} -
\Theta\bigr)^+
+ \mathbb{P}^{\pi} \bigl( \widetilde{\tau}_{\pi_{\infty}(c)} < \Theta\bigr)
\le c   \mathbb{E}^{\pi} (\tau- \Theta)^+
+ \mathbb{P}^{\pi} ( \tau< \Theta) .
\]
Since $\mathbb{P}^{\pi} ( \widetilde{\tau}_{[\pi_{\infty}(c),1]}
< \Theta) =
\alpha$, its expected detection delay has to be minimal compared to
other stopping time $\tau$'s for which $\mathbb{P}^{\pi} ( \tau<
\Theta) \le
\alpha$.

In this section, we show that $c \mapsto\mathbb{E}^{\pi} h ( \Pi_{
\widetilde
{\tau}_{\pi_{\infty}(c) } } )$ is a continuous function of $c \in
(0,\infty)$ with limits $0$ $(0)$ and $1-\pi$ $(1-p)$ as $c \searrow0$
and $ c \nearrow\infty$, respectively if $\pi> 0$ ($\pi=0$). Hence,
for a given pair $(\pi, \alpha)$ the arguments in \cite{Sh78} work,
and $\widetilde{\tau}_{\pi_{\infty}(c)}$ is optimal for the value of
$c$, for which $\mathbb{E}^{\pi} h ( \Pi_{ \widetilde{\tau}_{\pi
_{\infty}(c) }
} ) = \alpha$.

\subsection{False alarm probabilities} For a given threshold $r \in
(0,1)$, let $\widetilde{\tau}_r := \inf\{ t \ge0 \dvtx\Pi_t \ge r \}$ be
the exit time of $\Pi$ from the interval $[0,r)$,
and let
\begin{equation}
\label{defFr}
F_r ( \pi) := \mathbb{P}^{\pi} ( \widetilde{\tau}_{ r
} < \Theta)
= \mathbb{E}^{\pi} h ( \Pi_{ \widetilde{\tau}_{r } } )
\end{equation}
be the corresponding false alarm probability. On the event $ \{
\widetilde{\tau}_{r } < T_1 \}$, the exit time of $\Pi$ coincides
with the exit time $\tau_r$ of the process $Y$ [in \eqref{dynamicsofY}],
and we have
$ h ( \Pi_{ \widetilde{\tau}_{r } } ) = h ( Y_{
\tau_{r
} } ) $.
On the other hand, conditioned on $\{\widetilde{\tau}_{r } \ge T_1 \}$,
strong Markov property implies that the false alarm probability should
be computed by evaluating the function $F_r(\cdot)$ at the point $\Pi
_{T_1}$. Therefore, we expect the function $F_r(\cdot)$ to solve 
\begin{eqnarray*}
F_r (\pi) &=& \mathbb{E}^{\pi} \bigl[ 1_{ \{ \tau_{r } < T_1 \} } h
( Y_{
\tau_{r } } )
+ 1_{ \{ \tau_{r} \ge T_1 \} } \cdot F_r \bigl( Y_{T_1-} + p (1-
Y_{T_1-}) \bigr) \bigr]
\\
&= &\mathbb{E}^{\pi} \biggl[ e^{ - \lambda\tau_{r }} h ( Y_{
\tau_{r } }
)
+ \int_0^{\tau_r} \lambda e^{ - \lambda t} \cdot F_r(\mathbb
{S}(Y_t)) \,dt
\biggr] =: H^{(0)}_r [F_r] (\pi),
\end{eqnarray*}
where $H^{(0)}_r [\cdot] (\cdot)$ denotes the operator $H_r [\cdot]
(\cdot)$ in (\ref{defHr})--(\ref{Hrexplicit})
with $c=0$ [see also \eqref{defuandC}].
Hence,
if we apply the operator $H^{(0)}_r [\cdot] (\cdot)$ successively
starting with a suitably selected initial function, the sequence that
we obtain should convergence to the function $F_r(\cdot)$. Indeed, in
Appendix~\ref{secconstructingexittimeexpectation}, we show that the
sequence constructed as
\begin{equation}
\label{defus}
u_{0,r} (\cdot) = h(\cdot) \quad  \mbox{and}\quad
u_{n+1,r} (\cdot) = H^{(0)}_r [ u_{n, r}](\cdot)\qquad  \mbox{for $n
\in\mathbb{N}$,}
\end{equation}
is nonincreasing and converges uniformly to $F_r(\cdot)$ with error bounds
\begin{equation}
\label{uniformconvergenceofus}
0 \le F_r(\pi) \le u_{n ,r} (\pi) \le F_r(\pi) + (1-p)^n (1-\pi)\qquad
 \mbox{for all $n \in\mathbb{N}$.}
\end{equation}

It can easily be verified that the results in Lemmas \ref{lemmexittimeexpectation1} and
\ref{lemmexittimeexpectation2}
still hold for $c=0$. Hence, on the region $\{ (\pi, r ) \dvtx \pi< r \}$,
$u_{n,r} (\pi)$
has the form
\begin{eqnarray}
\label{unr}
&&\psi(\pi) \biggl( - \frac{ \eta(r)}{\psi(r)} \int_0^{r}
\frac{ u_{n-1,r}(\mathbb{S}(y)) }{ m_1-m_2 } \psi''(y) \,dy \nonumber\\
&&\hspace*{29pt}{}+
\int_{\pi}^{r} \frac{ u_{n-1,r}(\mathbb{S}(y)) }{ m_1-m_2 } \eta''(y)
\,dy + \frac{ h(r)}{\psi(r)}
\biggr)\\
&&\hspace*{7pt}\qquad{}+ \eta(\pi) \int_0^{\pi} \frac{ u_{n-1,r}(\mathbb
{S}(y))}{(m_1-m_2) }
  \psi''(y) \,dy ,\nonumber
\end{eqnarray}
thanks to identities $\mathcal{A}_0 \psi(\cdot) = \lambda\psi
(\cdot)$ and
$\mathcal{A}_0 \eta(\cdot) = \lambda\eta(\cdot)$. On the region
$\{ (\pi, r
) \dvtx\break \pi\ge r \}$, obviously, we have $u_{n,r} (\pi) = 1- \pi$.

\begin{lemm}
\label{lemmcontinuityofuinr}
For each $\pi\in[0,1]$, the functions $r \mapsto u_{n ,r} (\pi)$, for
$n \in\mathbb{N}$, and $r \mapsto F_{r} (\pi)$ are continuous on $r
\in(0,1)$.
\end{lemm}
\begin{pf}
The result is obvious for $\pi=1$, since $u_{n,r} (1)= F_r(1) = 1$,
for all $r \in(0,1)$. To prove the result for $\pi< 1$, we will show
that $(\pi,r) \mapsto u_{n,r}(\pi) $ is jointly continuous on $(0,1)
\times(0,1)$. However, observe that $u_{n,r} (\pi) = 1- \pi$, for
$\pi
\le r $; and $u_{n,r} (0) = u_{n-1,r} (p) $, for $r > \pi= 0 $ thanks
to \eqref{unr} [see also \eqref{continuityat0} in Lemma~\ref
{lemmexittimeexpectation2}]. Then direct computation gives
\[
\lim_{r \to0^+} u_{n,r} (0) = \lim_{r \to0^+} u_{n-1,r} (p) = 1- p <
1 = u_{n,0} (0),
\]
which shows that $(\pi,r) \mapsto u_{n,r}(\pi)$ is not continuous at $(0,0)$.

Clearly, $(\pi,r) \mapsto u_{0,r}(\pi)= 1- \pi$ is continuous on $(0,1)
\times(0,1)$. Suppose that the result holds for some $n \in\mathbb
{N}$. On
the region $\{ (\pi, r ) \dvtx \pi\ge r \}$, $u_{n+1,r}(\pi) $ again
equals $ h(\pi)$, and continuity is immediate. 

Also, using the joint continuity on $(0,1) \times(0,1)$ of the bounded
function $u_{n,r}(\pi) $
[and the boundary conditions $\psi'(0+) = 0$ and $\eta'(1-) = 0$],
it can be verified that the expression in \eqref{unr} is jointly
continuous on
$\{ (\pi, r ) \in(0,1) \times(0,1): \pi\le r \}$.
When we let $r \to\pi$ in \eqref{unr}, direct computation gives
\begin{eqnarray*}
&&\psi(\pi) \biggl( - \frac{ \eta(\pi)}{\psi(\pi)} \int_0^{\pi}
\frac{
u_{n,r}(\mathbb{S}(y))}{(m_1-m_2) }   \psi''(y)\, dy + \frac{
h(\pi)}{\psi(\pi)} \biggr) \\
&&\qquad{}+ \eta(\pi) \int_0^{\pi} \frac{
u_{n,r}(\mathbb{S}(y))}{(m_1-m_2) }   \psi''(y)\, dy = h(\pi).
\end{eqnarray*}
This implies that $u_{n+1,r}(\pi) $ is jointly continuous on $(0,1)
\times(0,1)$, and the result is true all $n \in\mathbb{N}$ by induction.

For $\pi=0 $ and $n \in\mathbb{N}$, we have $ u_{n+1,r}(0) =
u_{n,r}(p)$, and
the continuity of $r \mapsto u_{n+1,r}(0)$ follows from the first part
of the proof.
Finally, the uniform convergence in \eqref{uniformconvergenceofus}
imply that
$r \mapsto F_{r} (\pi)$ is also continuous, for each $\pi\in[0,1]$,
and this concludes the proof.
\end{pf}

By the definition of $F_{r}(\pi)$ given in \eqref{defFr}, we have
\begin{eqnarray}
\label{limitsofFraroundzero}
\lim_{r \to0^+} F_{r} (\pi) &=& 1- \pi\qquad
 \mbox{for } \pi> 0   \quad      \mbox{and}
 \nonumber
 \\[-8pt]
 \\[-8pt]
 \nonumber
 \lim_{r \to0^+} F_{r} (0) &=& 1- p ,
\end{eqnarray}
where the second limit follows from the behavior of the process $\Pi$
at $\{0\}$. That is, if $\Pi_0=0$, the process remains at this point
until the first arrival time $T_1$, and then it jumps to the point $\{
p\}$ [see \eqref{eqdynamicsofpi}].
Also note that, for all $\pi\in[0,1]$ and $r < 1$, the exit time
$\widetilde{\tau}_r$ is finite $\mathbb{P}^{\pi}$-almost surely,
and $\Pi
_{\widetilde{\tau}_r} \in(r , r + p (1-r))$. Hence,
\begin{equation}
\label{limitsofFraroundone}
\lim_{r \to1^-} F_{r} (\pi) = \lim_{r \to1^-} \mathbb{E}^{\pi}
[ 1- \Pi
_{\widetilde{\tau}_r}] = 0\qquad  
 \mbox{for $\pi\ge0$.}
\end{equation}

\begin{rem}
\label{remoptimalthresholdcontinuousinc}
The optimal threshold of the Bayesian formulation is a nonincreasing
and continuous function of the cost parameter $c$. If we let $\pi
_{\infty} (c)$ denote the optimal threshold as a function of $c$, we
have 
\begin{equation}
\label{limitsofpiinftyinc}
\lim_{c \to0^+ } \pi_{\infty}(c) = 1 
 \quad \mbox{and}\quad
\lim_{c \to\infty} \pi_{\infty}(c) = 0.
\end{equation}
\end{rem}

The limits in \eqref{limitsofpiinftyinc} can be obtained using the
bounds in \eqref{boundsonr}. Monotonicity of $\pi_{\infty}(c)$ in $c$
is also obvious and follows from \eqref{defR} and Remark~\ref
{remboundsonr}. For the proof of the continuity of $c \mapsto\pi
_{\infty} (c)$,
Appendix~\ref{secappendixB} can be consulted.

Lemma~\ref{lemmcontinuityofuinr} and Remark~\ref
{remoptimalthresholdcontinuousinc} imply that 
$F_{\pi_{\infty}(c)} (\pi)$ 
is continuous with respect to $c$ on $(0, \infty)$. Moreover, thanks to
(\ref{limitsofFraroundzero})--(\ref{limitsofFraroundone}) we have
\begin{eqnarray*}
\lim_{c \to0^+} F_{\pi_{\infty}(c)} (\pi) &=& \lim_{r \to1-} F_{r}
(\pi) = 0 ,\qquad  \mbox{with} \\
\lim_{c \to\infty} F_{\pi_{\infty}(c)} (\pi) &=&
\lim_{r \to0+} F_{r} (\pi) = 1- \pi 
\end{eqnarray*}
for $\pi> 0$, and
\begin{eqnarray*}
\lim_{c \to0^+} F_{\pi_{\infty}(c)} (0) &=& \lim_{r \to1-} F_{r}
(0) = 0 ,\qquad  \mbox{with} \\
\lim_{c \to\infty} F_{\pi_{\infty}(c)} (0) &=& \lim_{r \to0+}
F_{r} (0) = 1- p .
\end{eqnarray*}
Hence (excluding the trivial cases) it is possible to pick a value of
$c$ such that the exit time $\widetilde{\tau}_{\pi_{\infty}(c)}$
has a
false alarm probability $\alpha$ and solves the variational formulation.

\begin{appendix}

\section{On the conditional probability process}
\label{secappendixA}

\subsection{\texorpdfstring{An auxiliary probability measure and the proof of \protect\eqref
{eqPiandPhi}}{An auxiliary probability measure and the proof of (2.7)}} Let $(\Omega, \mathcal{H},\break \mathbb{P}_0)$ be a
probability space
hosting the following independent stochastic elements:
\begin{itemize}
\item[$\bullet$] a Wiener process $X$ (with $\mu=0$),
\item[$\bullet$] a simple Poisson process $N$ with arrival rate $\lambda$ and
arrival times $(T_n)_{n \ge0}$,
%
\item[$\bullet$] an integer valued random variable with distribution
$ \mathbb{P}_0 \lbrace\zeta= 0 \rbrace= \pi$ and
$ \mathbb{P}_0 \lbrace\zeta= n \rbrace=(1- \pi)
(1-p)^{n-1}p $
for $n \in\mathbb{N}$,
\item[$\bullet$] a random variable $\Theta$ defined as in \eqref{defTheta}.
\end{itemize}

Let $\mathbb{G} \equiv \{ \mathcal{G} \}_{t \ge0}$ be an extended
filtration such that
$\mathcal{G}_t := \sigma\{ X_s, N_s , \zeta\dvtx s\leq t\} $. In terms of the
process $L_t = \exp\{ \mu X_t - \mu^2 t / 2\}$, we introduce a new
probability measure $\mathbb{P}$ whose Radon--Nykodyn derivative is
\[
Z_t := \frac{d \mathbb{P}}{ d \mathbb{P}_0 }\bigg|_{\mathcal{G}_t}
= 1_{ \{ \Theta>
t \} }
+ 1_{ \{ \Theta\le t \} } \frac{L_t}{L_{\Theta}}.
\]
Under the new measure, the process $X$ is a Brownian motion that gains
a drift~$\mu$ at~$\Theta$. The random variables $\zeta$ and $\Theta$
have the same distribution under $\mathbb{P}$ since $\zeta\in
\mathcal{G}_0$
and $Z_0 =1$. In other words, we have the same setup described in
Sections~\ref{secintro} and~\ref{secproblem}.

Let us now define the likelihood ratio process
\[
\Phi_t := \frac{ \mathbb{P}\{ \Theta\le t | \mathcal{F}_t \} }{
\mathbb{P}\{ \Theta> t | \mathcal{F}
_t \}}
= \frac{\mathbb{E}_0 [ Z_t 1_{\{ \Theta\le t\} } | \mathcal{F}_t ] }
{\mathbb{E}_0 [ Z_t 1_{\{ \Theta> t\} } | \mathcal{F}_t ]},
\]
where the equality follows from Bayes' rule. Using the independence
of $X,N$ and $\zeta$ under $\mathbb{P}_0$, we obtain
\[
\mathbb{E}_0 \bigl[ Z_t 1_{\{ \Theta\le t\} } | \mathcal{F}_t \bigr] =
\pi L_t + (1-\pi) \sum_{i=1}^{N_t} (1-p)^{i-1} p \frac{L_t}{L_{T_i}}
\]
and
\[
\mathbb{E}_0 \bigl[ Z_t 1_{\{ \Theta> t\} } | \mathcal{F}_t \bigr] =
\mathbb{P}_0 \bigl[ 1_{\{ \Theta> t\} } | \mathcal{F}_t \bigr] = (1-\pi) (1-p)^{N_t}.
\]
Therefore, we have
\[
\Phi_t = \frac{L_t}{(1-p)^{N_t}} \Biggl( \frac{\pi}{1-\pi} +
\sum_{i=1}^{N_t} \frac{(1-p)^{i-1} p }{L_{T_i}} \Biggr),
\]
and this proves \eqref{eqPiandPhi}.

\subsection{Constructing the exit time (false alarm) probabilities}
\label{secconstructingexittimeexpectation}
Let $H^{(0)}_{r}$ denotes $H_{r}$ defined in \eqref{defHr} with
$c=0$. It should be noted that the proofs of Lemmas~\ref
{lemmexittimeexpectation1} and \ref{lemmexittimeexpectation2}
use only the continuity of the given function $w(\cdot)$ and the bounds
$0 \le w(\cdot) \le h(\cdot)$. Hence, they also cover the case $c=0$.

\begin{rem}
\label{rempropertiesofHzeror}
The operator $H^{(0)}_{r}$ is monotone in $w(\cdot)$; that is for $w_1
(\cdot) \le w_2(\cdot) $, we have $H^{(0)}_{r}[w_1](\cdot) \le
H^{(0)}_{r}[w_2](\cdot)$. Moreover, if $w(\cdot)$ is a continuous
function bounded as $0 \le w(\cdot) \le h(\cdot)$, then so is
$H^{(0)}_{r} [w](\cdot)$.
\end{rem}
\begin{pf}
The claim on monotonicity is obvious. Given $w(\cdot) $ continuous and
bounded as $0 \le w(\cdot) \le h(\cdot)$, $H^{(0)}_{r} [w](\cdot)$ is
again continuous by Lemma~\ref{lemmexittimeexpectation2}.

Since the process $Y$ in \eqref{dynamicsofY} 
is a bounded martingale, we have
\begin{eqnarray*}
\ell(\pi) &:=& \mathbb{E}^{\pi} \int_0^{\infty} e^{- \lambda t }
\lambda
h (\mathbb{S}(Y_t)) \,dt\\
&=&
\mathbb{E}^{\pi} \int_0^{\infty} e^{- \lambda t } \lambda(1-p) (1-
Y_t)\, dt\\
& =&
(1-p) (1- \pi).
\end{eqnarray*}
Then, for a function $w(\cdot)$ bounded as $0 \le w(\cdot) \le
h(\cdot) $, strong Markov property gives
\begin{eqnarray*}
0 &\le& H^{(0)}_{r}[w](\pi)
\le\mathbb{E}^{\pi} \biggl[ e^{- \lambda\tau_r } h (Y_{\tau_r} )
+ \int
_0^{\tau_r } e^{- \lambda t } \lambda  h (\mathbb{S}(Y_t)) \,dt
\biggr]
\\
&=& \ell(\pi) + \mathbb{E}^{\pi} e^{- \lambda\tau_r } [ h
(Y_{\tau_r} )
- \ell(Y_{\tau_r} ) ] \\
&=& \ell(\pi) + \mathbb{E}^{\pi} e^{- \lambda\tau_r } p (1 -
Y_{\tau_r} ) \le
\ell(\pi) + \mathbb{E}^{\pi} p (1 - Y_{\tau_r} ) = h (\cdot).
\end{eqnarray*}
Hence, $0 \le H^{(0)}_{r}[w](\cdot) \le h(\cdot) $ again.
\end{pf}

Using Remark~\ref{rempropertiesofHzeror} above, it can be shown by
induction (as in the proof of Remark~\ref{vnsconverge}) that the sequence
\begin{eqnarray}
\label{sequencesus}
u_{0,r} (\cdot) = h(\cdot) \quad  \mbox{and}\quad
u_{n+1,r} (\cdot) = H^{(0)}_r [ u_{n, r}](\cdot)\qquad  \mbox{for $n
\in\mathbb{N}$,}
\end{eqnarray}
is nonincreasing, 
and each function is nonnegative, continuous and bounded above by
$h(\cdot) $. The pointwise limit $u_{\infty,r} (\cdot) := \inf_{n
\in
\mathbb{N}} u_{n,r} (\cdot)$ exists and it is bounded as $0 \le
u_{\infty,r}
(\cdot) \le h(\cdot)$.

\begin{rem}
The limit function $u_{\infty,r} (\cdot) $ solves
$ u_{\infty,r} (\cdot) = H^{(0)}_r [ u_{\infty,r}] (\cdot)$, on $[0,1]$.
\end{rem}
\begin{pf}
The proof follows from a straightforward modification of \eqref
{fixedpointofJ1} by replacing $v_{\infty} $, $v_n $, $\tau$ with
$u_{\infty,r} $, $u_{n,r} $, $\tau_r$ respectively.
\end{pf}

\begin{rem}
\label{remuniformconvergenceofus}
The sequence defined in \eqref{sequencesus} converges uniformly on $
[0,1]$, and we have the explicit error bounds
\begin{equation}
\label{uniformconvergenceofurs}
0 \le u_{n, r }(\pi) - u_{\infty,r} (\pi) \le(1-p)^n (1- \pi)
 \qquad \mbox{for $n \in\mathbb{N}$.}
\end{equation}
\end{rem}
\begin{pf}
We will establish the inequalities above by modifying the proof of
Lemma~\ref{lemmuniformconvergence}.

The first inequality in \eqref{uniformconvergenceofus} is obvious.
The second inequality follows immediately for $n=0$ since $0 \le
u_{\infty,r} (\cdot) \le h (\cdot) $. Assume it holds for some $n
\in\mathbb{N}
$. Then using the induction hypothesis and the identity $ u_{\infty,r}
(\cdot) = H^{(0)}_r [ u_{\infty,r}] (\cdot)$, we have
\begin{eqnarray*}
u_{n+1,r}
(\pi)&=& H^{(0)}_r [ u_{n, r}](\pi) \\
&\le&\mathbb{E}^{\pi} \biggl[ e^{- \lambda\tau_r } h (Y_{\tau_r} ) \\
&&\hspace*{17pt}{}+
\int_0^{\tau_r
} e^{- \lambda t } \lambda[ u_{\infty, r} (\mathbb{S}(Y_t)) +
(1- p
)^{n+1} (1- Y_t)]\, dt \biggr] \\
&\le& u_{\infty,r} (\cdot) + \mathbb{E}^{\pi} \biggl[ \int
_0^{\infty} e^{-
\lambda t } \lambda[ (1- p )^{n+1} (1- Y_t)] \,dt
\biggr] \\
&=& u_{\infty,r} (\cdot) + (1- p )^{n+1} (1- \pi),
\end{eqnarray*}
%
and \eqref{uniformconvergenceofurs} follows.
\end{pf}

\begin{cor}
Since, each $u_{n,r} (\cdot)$ is continuous, so is $u_{\infty,r}
(\cdot)$ thanks to Remark~\textup{\ref{remuniformconvergenceofus}}. Then,
the identity $ u_{\infty,r} (\cdot) = H^{(0)}_r [ u_{\infty,r}]
(\cdot
) $ and Lemma~\textup{\ref{lemmexittimeexpectation2}} imply that the
function $u_{\infty,r} (\cdot)$ solves
\begin{equation}
( - \lambda+ \mathcal{A}_0 )u_{\infty,r}(\pi) + \lambda u_{\infty
,r}(\mathbb{S}(y)) =0\qquad
 \mbox{on $(0,r)$,}
\end{equation}
and at $\pi=0$, we have $u_{\infty,r}(0)= u_{\infty,r}(p)$ [see
\eqref{continuityat0}].
\end{cor}

\begin{prop}
\label{propuinftyequalsFr}
The limit function $u_{\infty,r}(\cdot)$ coincides on $[0,1]$ with the
exit time expectation $F_r(\cdot)$ defined in \eqref{defFr}.
\end{prop}
\begin{pf}
The characterization in \eqref{Hrexplicit} indicates that the
derivative of $u_{\infty,r}$ is bounded on $(l,r)$, for $0< l < r$.
Then, for $\pi\in(l,r)$, a localization argument and It\^{o}'s rule gives
\begin{eqnarray}
\label{Itoforuinfty}
&&\mathbb{E}^{\pi} u_{\infty,r} \bigl(\Pi_{\widetilde{\tau}_{[l,r]}} \bigr)
 \nonumber\\
&&\qquad=u_{\infty,r}(\pi) + \mathbb{E}^{\pi} \int_0^{ \widetilde{\tau
}_{[l,r ]} }
[ (- \lambda+ \mathcal{A}_0) u_{\infty,r} (\Pi_{u-})
+ \lambda u_{\infty,r} (\mathbb{S}(\Pi_{u-})) ]\, du \\
&&\qquad= u_{\infty
,r}(\pi),\nonumber
\end{eqnarray}
where $\widetilde{\tau}_{[l,r]}$ is the exit time of $\Pi$ from the
interval $(l,r)$. The boundary $\{ 0 \}$ is natural for the
diffusive part of the process $\Pi$ and its jumps are positive
(toward $\{ 1 \}$). This implies that $\widetilde{\tau}_{[l,r]}
\nearrow\widetilde{\tau}_{r} = \inf\{ t \ge0 \dvtx \Pi_t \ge r \}$ as
$l \to0^+$, $\mathbb{P}^{\pi}$-almost surely [see also~\eqref{exittimeexpectationuniformlybounded}]. Therefore, when we
let $l \to0^+$ in \eqref{Itoforuinfty} we obtain
\begin{eqnarray*}
u_{\infty,r}(\pi) &=& \lim_{l \to0^+} \mathbb{E}^{\pi} u_{\infty
,r} \bigl(\Pi
_{\widetilde{\tau}_{[l,r]}} \bigr)
\\
&=& \lim_{l \to0^+} u_{\infty,r} ( l )
\mathbb{P}^{\pi} \bigl\{ \widetilde{\tau}_{[l,r]} < \widetilde{\tau
}_{r} \bigr\}
+ \mathbb{E}^{\pi} 1_{ \{ \widetilde{\tau}_{[l,r]} = \widetilde
{\tau}_{r} \}} h
( \Pi_{\widetilde{\tau}_{r} } )\\
& =& \mathbb{E}^{\pi} h (\Pi
_{\widetilde{\tau}_{ r
}} ).
\end{eqnarray*}
This shows $u_{\infty,r}(\cdot) = F_r (\cdot)$ on $(0,r)$.

When $\Pi_0 = 0 $, the process stays at $\{ 0 \}$ until the first
arrival time $T_1$ of $N$. It jumps to $\{ p \}$ at $T_1$. Hence, by
strong Markov property, we have $F_r (0) = F_r (p)$, and this shows
$u_{\infty,r}(0) = F_r (0)$ [since $u_{\infty,r}(0) = u_{\infty
,r}(p)$]. Finally, for $\pi\ge r $, we have $u_{\infty,r}(\pi) = 1-
\pi
$ by the construction in \eqref{sequencesus}; hence, the equality
$u_{\infty,r}(\cdot) = F_r (\cdot)$ is obvious.
\end{pf}

\section{Other proofs}
\label{secappendixB}

\begin{pf*}{Proof of \protect\eqref{eqhatWandNindependent}}
The process $\widehat{W}$ is a
$(\mathbb{P}, \mathbb{F})$-Brownian motion (this can be verified
using L\'{e}vy's
characterization for Brownian motion) and $N$ is a $(\mathbb{P},
\mathbb{F})$-Poisson
process. Therefore, it is sufficient to show \eqref
{eqhatWandNindependent} for $s_1 = s_2 = s$.

Note that the process $\widehat{W}$ can be written as
\begin{equation}
\label{barWdecomposed}
\widehat{W}_t = W_t + \mu\int_0^t \bigl[ 1_{ \{ \Theta\le u \} } -
\Pi
_u \bigr] \,du .
\end{equation}
Therefore, if we apply It\^{o} formula to real and imaginary parts of
the process
$K_t := f(\widehat{W}_t , N_t)$, for
$
f(x,y) = \exp\{ir x + i qy \}
$,
we obtain
\begin{eqnarray}
\label{ItoforK}
K_s &=& K_t
+ i \biggl[ \int_t^s r K_u \,d W_u + \int_t^s r \mu K_u 1_{ \{ \Theta
\le u \} }\, d u - \int_t^s r \mu K_u \Pi_u\, d u \biggr]\nonumber\\
&&{}- \frac{1}{2} r^2 \int_t^s K_u\, du
+ \int_t^s (e^{iq}-1) K_u ( d N_u - \lambda\, du )\\
&&{}+ \int_t^s \lambda(e^{iq}-1) K_u\, du\nonumber
\end{eqnarray}
for $t \le u \le s$. Clearly, we have
\[
\mathbb{E}\biggl[ \int_t^s (e^{iq}-1) K_u ( d N_u - \lambda\, du )
\Big| \mathcal{F}
_t \biggr] = 0 =
\mathbb{E}\biggl[ \int_t^s K_u \,d W_u   \Big| \mathcal{F}_t \biggr].
\]
Moreover, for a set $A \in\mathcal{F}_t $ we have $\mathbb{E}1_A K_u
1_{ \{ \Theta
\le u \} } = \mathbb{E}1_A K_u \Pi_u$.\vspace*{1.5pt} Then by multiplying both sides in
\eqref{ItoforK} with $1_A / K_t = 1_A \cdot e^{- i r \widehat{W}_t
- i q N_t }$ and taking the expectations we get
\begin{eqnarray*}
&&\mathbb{E}[ 1_A \exp\{ i r( \widehat{W}_s - \widehat
{W}_t) + i q (
N_s - N_t) \} ]
\\
&&\qquad= P(A) + \int_t^s \biggl( - \frac{r^2}{2}+\lambda(e^{iq}-1)
\biggr)\\
&&\qquad\hspace*{63pt}{}\times
\mathbb{E}[ 1_A \exp\{ i r( \widehat{W}_u - \widehat
{W}_t )
+ i q ( N_u - N_t) \} ]\, du .
\end{eqnarray*}
By solving this integral equation for the (deterministic) function
\[
\varrho_t (\cdot) \dvtx s \mapsto\mathbb{E}[ 1_A \exp\{ i r(
\widehat{W}_s - \widehat{W}_t) + i q ( N_s - N_t)\} ]
\]
we obtain
$\varrho_t (s) =P(A) \cdot\exp\{ ( - \frac
{r^2}{2}+\lambda(e^{iq}-1)
) (s-t) \}$,
and this proves \eqref{eqhatWandNindependent} for $s_1=s_2 =s$.
\end{pf*}

\begin{pf*}{Proof of \protect\eqref{inequalityonstrictconcavity}}
Let $\pi$ be a fixed point on $(0,r[w])$. For any $r \ge r[w]$,
$H_r[w](\pi)$ is given by
\eqref{Hrexplicit} with
\begin{eqnarray*}
\frac{ \partial H_r[w](\pi)}{\partial r}
&=& - \frac{\psi(\pi) }{\psi^2(r)} \biggl\lbrace h(r) \psi'(r) +
\psi(r) - (m_1 - m_2 ) \int_0^r u_2[w](y) \,dy \biggr\rbrace\\
&=& - \frac{\psi(\pi) }{\psi^2(r)} B[w](r).
\end{eqnarray*}
The last expression is strictly positive for $r > r[w]$ since $B[w](r)$
is strictly negative thanks to Lemma~\ref{lemmuniqueroot}. This
implies that $H_{r[w]}[w](\pi) < H_{r_1}(\pi) < H_{r_2}(\pi)$, for all
$r[w] < r_1< r_2$, and we have $H_r[w](\pi) < \lim_{r \nearrow1}
H_r[w](\pi)$. Since the right boundary is natural, $\tau_r \nearrow
\infty$ as $r \nearrow1$. Then by dominated convergence theorem (see
Remark~\ref{remintegrabilitycond}), we obtain
\begin{eqnarray*}
H_r[w](\pi) &<& \mathbb{E}^{\pi} \biggl[ \int_0^{\infty} e^{-
\lambda t } \bigl(
g(Y_t ) + \lambda w(\mathbb{S}(Y_t))\bigr)\, dt
\biggr] \\
&\le&\int_0^{\infty} \bigl( g(\pi) + \lambda w(\mathbb{S}(\pi))
\bigr)\, dt\\
&=& \frac{g(\pi) + \lambda w(\mathbb{S}(\pi))}{\lambda},
\end{eqnarray*}
where the second inequality is by Jensen's inequality [recall that
$Y$ is a martingale and $w(\cdot)$ is concave], and
\eqref{inequalityonstrictconcavity} follows.
\end{pf*}

\begin{pf*}{Proof of Remark~\protect\ref{remoptimalthresholdcontinuousinc}}
The limits in \eqref{limitsofpiinftyinc} follow easily from
\eqref
{boundsonr}. It is also clear that $c \mapsto\pi_{\infty}(c)$ is
nonincreasing thanks to \eqref{defR} and Remark~\ref
{remboundsonr}. Here, we show that $c \mapsto\pi_{\infty}(c)$ is
continuous on $(0,\infty)$.

Let $V_c(\pi)$ and $B_c[\cdot]$ denote respectively the dependence on
$c$ of the value function $V$ and the operator $B[\cdot]$ defined in
\eqref{defB}.

Since the value function $V$ is a fixed point of the operator $J$,
Lemma~\ref{lemmuniqueroot} gives
\[
B_{c_1}[V_{c_1}] (\pi_{\infty}(c_1) ) =0 = B_{c_1}[V_{c_1}] (\pi
_{\infty}(c_2) )\qquad\mbox{for $0 < c_1 \le c_2 < \infty$.}
\]
By using these equalities
together with the explicit form of $B_{\cdot}[\cdot]$ in \eqref{defB}
[and the identity $\mathcal{A}_0 \psi(\cdot)= \lambda\psi(\cdot)$],
we obtain
\begin{eqnarray}
\label{matchingBs}
0 &\le& B_{c_1}[V_{c_1}] (\pi_{\infty}(c_2) ) - B_{c_1}[ V_{c_1} ]
(\pi_{\infty}(c_1) ) \nonumber\\
&\equiv&\int_{ \pi_{\infty}(c_2) }^{\pi_{\infty}(c_1) }
\frac{\psi''(y)}{\lambda} [ c_1 y + \lambda V_{c_1}(\mathbb{S}(y))
- \lambda h(y) ] \,dy \\
&=& \int_0^{\pi_{\infty}(c_2) } \frac{\psi''(y)}{\lambda}
\bigl[ (c_2 - c_1) y + \lambda\bigl( V_{c_2}(\mathbb{S}(y)) - V_{c_1}
(\mathbb{S}(y)) \bigr)
\bigr]\, dy.\nonumber
\end{eqnarray}
%
%
Moreover, we have
$V_{c_2}(\pi) \le
\mathbb{P}^{\pi} ( \widetilde{\tau}_{ \pi_{\infty}(c_1) } <
\Theta) +
{c_2} \mathbb{E}^{\pi} ( \widetilde{\tau}_{ \pi_{\infty}(c_1) } -
\Theta)^+
$, and this gives
the Lipschitz condition
\[
\frac{V_{c_2}(\pi) - V_{c_1}(\pi)}{ c_2-c_1 }
\le
\mathbb{E}^{\pi} \bigl( \widetilde{\tau}_{ \pi_{\infty}(c_1) } -
\Theta\bigr)^+
\le\frac{ V_{c_1}(\pi) }{ c_1 } \le\frac{1}{\delta}
\]
for any $\delta< c_1 $.
Using this inequality in \eqref{matchingBs}, we obtain
\begin{eqnarray}
\label{boundwithc1c2}
0 &\le&\int_{ \pi_{\infty}(c_2) }^{\pi_{\infty}(c_1) } \frac{\psi
''(y)}{\lambda}
[ c_1 y + \lambda V_{c_1}(\mathbb{S}(y)) - \lambda h(y) ]\, dy
\nonumber\\
&\le&
\int_0^{\pi_{\infty}(c_2) } \frac{\psi''(y)}{\lambda}
(c_2 - c_1) \biggl[ y + \frac{\lambda}{\delta} \biggr]\, dy\\
&\le&(c_2 - c_1) \biggl[ \frac{1}{\lambda} + \frac{1}{\delta}
\biggr] \psi
'( \pi_{\infty}(c_2) ) .\nonumber
\end{eqnarray}
This implies that $ \pi_{\infty}(c_2) \nearrow\pi_{\infty}(c_1)$ as
$c_2 \searrow c_1$.

Similarly, it is easy to show that
\begin{eqnarray*}
0 &\le&\int_{ \pi_{\infty}(c_2) }^{\pi_{\infty}(c_1) }
\frac{\psi''(y)}{\lambda} [ c_2 y + \lambda V_{c_2}(\mathbb{S}(y))
- \lambda h(y) ]\, dy \\
&=& \int_0^{\pi_{\infty}(c_1) } \frac{\psi''(y)}{\lambda}
\bigl[ (c_2 - c_1) y + \lambda\bigl( V_{c_2}(\mathbb{S}(y)) - V_{c_1}
(\mathbb{S}(y)) \bigr)
\bigr] \,dy\\
&\le&\int_0^{\pi_{\infty}(c_1) } \frac{\psi''(y)}{\lambda}
(c_2 - c_1) \biggl[ 1 + \frac{\lambda}{\delta} \biggr]\, dy\\
&=& (c_2 - c_1) \biggl[ \frac{1}{\lambda} + \frac{1}{\delta}\biggr]
\psi'(\pi_{\infty}(c_1))
\end{eqnarray*}
for some $\delta< c_1$.
This shows that $ \pi_{\infty}(c_1) \searrow
\pi_{\infty}(c_2)$ as $c_1 \nearrow c_2$, and the continuity of $c
\mapsto\pi_{\infty}(c)$ follows.
\end{pf*}
\end{appendix}

\section*{Acknowledgments}
The author is thankful to the referees and the associate editors for
their valuable remarks and suggestions, which improved the presentation
of this paper.

%

\printaddresses

\end{document}